\documentclass[10pt]{amsart}
\usepackage{amsmath, amsfonts, amssymb, amstext, amscd, amsthm}
\allowdisplaybreaks

\newtheorem{goins}{Theorem}[subsection]

\newtheorem{infinitesimal}[goins]{Proposition}
\newtheorem{goins_lemma}[goins]{Proposition}
\newtheorem{existence}[goins]{Proposition}
\newtheorem{dimension}[goins]{Proposition}

\newtheorem{localization}[goins]{Proposition}
\newtheorem{isomorphism}[goins]{Proposition}
\newtheorem{wiles}[goins]{Theorem}
\newtheorem{buzzard}[goins]{Theorem}

\newtheorem{klein}[goins]{Proposition}
\newtheorem{q-curve}[goins]{Proposition}
\newtheorem{ell_rep}[goins]{Proposition}
\newtheorem{modularity}[goins]{Proposition}

\begin{document}


\title{On the Modularity of Wildly Ramified Galois Representations}

\author{Edray Herber Goins}
\address{California Institute of Technology \\ Mathematics 253-37 \\ Pasadena, CA 91125}
\email{goins@caltech.edu}



\begin{abstract}
We show that an infinite family of odd complex 2-dimensional Galois representations ramified at 5 having nonsolvable projective image are modular, thereby verifying Artin's conjecture for a new case of examples.  Such a family contains the original example studied by Buhler.  In the process, we prove that an infinite family of residually modular Galois representations are modular by studying $\Lambda$-adic Hecke algebras.  \end{abstract}

\maketitle


\section{Introduction}

There is considerable interest in continuous homomorphisms
\[ \rho: G_{\mathbb Q} \to GL_2 \left( \overline {\mathbb Q}_\ell \right) \]

\noindent where $G_{\mathbb Q} = \text{Gal} \left( \overline {\mathbb Q} / \mathbb Q \right)$ is the absolute Galois group of $\mathbb Q$ and $\ell$ is a fixed rational prime.  For example, $\rho = \rho_{E, \ell}$ may be the $\ell$-adic representation of an elliptic curve $E$ over $\mathbb Q$, or $\rho = \rho_f$ may be the $\ell$-adic representation associated to a modular form.  The continuity of such Galois representations implies the image lies in $GL_2(\mathcal O)$ for some ring of integers $\mathcal O$ with maximal ideal $\lambda$ in a finite extension $K$ of $\mathbb Q_\ell$; then $k = \mathcal O / \lambda$ is a finite extension of $\mathbb F_\ell$.  We define the residual representation $\overline \rho$ as the composition
\[ \overline \rho: G_{\mathbb Q} \to GL_2 \left(\mathcal O \right) \to GL_2 \left( \overline {\mathbb F}_\ell \right). \]

\noindent  For example, $\overline \rho = \overline \rho_{E, \ell}$ may be the mod $\ell$ representation of an elliptic curve, or $\overline \rho = \overline \rho_f$ may be the mod $\ell$ reduction of the representation associated to a modular form.  We also consider the projective representation $\widetilde \rho$ as the composition
\[ \widetilde \rho : G_{\mathbb Q} \to GL_2 \left( \overline {\mathbb Q}_\ell \right) \to PGL_2 \left( \overline {\mathbb Q}_\ell \right). \]

A common question is: Given a continuous $\ell$-adic Galois representation $\rho$ such that $\overline \rho$ is modular i.e. $\overline \rho \simeq \overline \rho_f$, when is $\rho$ modular i.e. $\rho \simeq \rho_f$?  In this paper the main result is the

\begin{wiles} For $\ell$ an odd prime, let $\rho: G_{\mathbb Q} \to GL_2(\mathcal O)$ be a continuous $\ell$-adic representation such that
\begin{enumerate}
\item $\rho$ is ordinary and ramified at finitely many primes;
\item $\overline \rho$ is absolutely irreducible when restricted to $\text{Gal} \left( \overline {\mathbb Q} / \mathbb Q(\sqrt{(-1)^{(\ell-1)/2} \, \ell}) \right)$, modular, and wildly ramified at $\ell$.
\end{enumerate}

\noindent Then $\rho$ is $\ell$-adically modular i.e. $\rho \simeq \rho_f$ for an $\ell$-adic cusp form $f$.
\end{wiles}

The proof the theorem follows the ideas outlined in \cite{MR98c:11047} and \cite{MR1639612}, but we generalize to $\Lambda$-adic modular forms.  One application of this result is the

\begin{buzzard} For $\ell$ an odd prime, let $\rho: G_{\mathbb Q} \to GL_2 \left( \mathcal O \right)$ be a continuous Galois representation such that
\begin{enumerate}
\item $\rho$ is ramified at finitely many primes;
\item $\overline \rho$ is absolutely irreducible when restricted to $\text{Gal} \left( \overline {\mathbb Q} / \mathbb Q(\sqrt{(-1)^{(\ell-1)/2} \, \ell}) \right)$, modular, and wildly ramified at $\ell$;
\item $\rho(G_\ell)$ is finite and $\widetilde \rho(G_\ell)$ is a cyclic group of $\ell$-power order.
\end{enumerate}

\noindent Then $\imath \circ  \rho: G_{\mathbb Q} \to GL_2 \left( \mathbb C \right)$ is modular for each embedding $\imath: K \hookrightarrow \mathbb C$.  \end{buzzard}

The proof of this result uses new ideas outlined in \cite{MR1937198} on ``glueing'' together two $\ell$-adic modular forms.  Our paper is different from others in the literature (most notably \cite{MR2000j:11062}) in that we consider deformations of wildly ramified mod $\ell$ representations, not tamely ramified ones.  We use a trick that exploits the wild ramification of the representations at hand.

Our third result is the

\begin{goins} Let $\rho: G_{\mathbb Q} \to GL_2(\mathbb C)$ be a continuous representation  with nonsolvable image, and denote $L/\mathbb Q$ as the extension cut out by $\widetilde \rho$.  If this extension is the splitting field of $x^5 + B \, x + C$ such that $75 \, C^2 / \sqrt{256 \, B^5 + 3125 \, C^4}$ is the square of a 5-adic unit, then $\rho$ is (classically) modular.  \end{goins}

The proof relies on the modularity of $\mathbb Q$-curves, as shown in \cite{MR1844206} using ideas from \cite{MR1639612}.  One example that is covered in this result is the original example in \cite{MR58:22019}; the polynomials
\[ x^5 + 10 \, x^3 + 10 \, x^2 + 35 \, x + 18 \qquad \text{and} \qquad 5 \, x^5 + 20 \, x + 16 \]

\noindent generate the same splitting field.  An infinite family of examples are
\[ x^5 + 5 \left( \frac {9-5 \, u^4}{5 \, u^4} \right) x + 4 \left( \frac {9-5 \, u^4}{5 \, u^4} \right), \qquad u \in \mathbb Q \cap \mathbb Z_5^{\times}. \]

\noindent (When $u=1$ we find the previous example.)  This result gives the first known proof of infinitely many examples of icosahedral Galois representations satisfying Artin's Conjecture which are ramified at 5 since the conditions in the theorem force $\rho$ to be wildly ramified there.  Other papers in the literature consider representations which are unramified at 5: \cite{MR96a:11128}, \cite{MR2002c:11057}, and \cite{MR1897901} consider finitely many representations ramified at either 2 or 3 yet unramified at 5, while \cite{MR1845181} considers infinitely many representations unramified at both 2 and 5.  Only \cite{MR58:22019} considers a representation ramified at 5.

We organize this exposition as follows.  In section 2 we review the Galois cohomology involved in computing the dimension of the universal deformation ring.  We consider ordinary representations such that the restriction $G_\ell \to GL_2(\mathbb F_\ell)$ may have equal eigenvalues, so we explain why the universal deformation ring exists.  In section 3 we review the theory of ordinary $\Lambda$-adic modular forms and how this gives a family of modular Galois representations $G_{\mathbb Q} \to GL_2(\overline {\mathbb Q}_\ell [[X]])$, following closely \cite[Chap. 7]{MR94j:11044}.  We then show that one may prove certain weight 1 representations are modular if one first proves certain weight 2 representations are modular, and we employ standard ``lifting'' arguments as in \cite{MR96d:11071}, \cite{MR98c:11047} and \cite{MR1639612} to show the universal deformation ring is indeed modular in the weight 2 case.  In section 4 we present applications of the theorems above using ideas from \cite{MR18:329c} relating icosahedral Galois extensions with elliptic curves.

The author would like to thank Richard Taylor for invaluable help with this project.  The author would also like to thank Harvard University for support, as well as Kevin Buzzard, Mark Dickinson, Matthias Flach, and Chris Skinner for helpful conversations.

\tableofcontents


\section{Universal Deformation Ring}

Every $\ell$-adic Galois representation $\rho: G_{\mathbb Q} \to GL_2(\mathcal O)$ has a mod $\ell$ reduction $\overline \rho: G_{\mathbb Q} \to GL_2(k)$.  Conversely for a local $\mathcal O$-algebra $\mathcal O'$ which is complete, Noetherian, and has residue field $k' = \mathcal O' / \lambda'$, we say a representation $\rho': G_{\mathbb Q} \to GL_2(\mathcal O')$ is a deformation of $\overline \rho$ if $\overline \rho' \simeq \overline \rho \otimes_k k'$ and $\det \rho' = \det \rho$.  When $\rho$ ramifies at finitely many primes and $\overline \rho$ is absolutely irreducible, there exists a universal deformation $\rho^{\text{univ}}: G_{\mathbb Q} \to GL_2(R)$ i.e. for any deformation $\rho'$ there exists a unique $\mathcal O$-algebra surjection $\phi: R \to \mathcal O'$ such that $\rho' \simeq \phi \circ \rho^{\text{univ}}$; for more properties, see \cite{ MR90k:11057}.  Unfortunately, the universal deformation ring $R$ is too big to work with but we may consider a filtration of finitely generated deformation rings $R_\Sigma$ corresponding to finite sets $\Sigma$.  Our goal in this section is to give an explicit formula for the topological dimension of $R_\Sigma$ as an $\mathcal O$-algebra.
\subsection{Galois Theory}

We briefly review some preliminaries to fix notation.

Given a place $\nu$ of $\mathbb Q$, each embedding $\imath_\nu: \overline {\mathbb Q} \hookrightarrow \overline {\mathbb Q}_\nu$ induces an embedding $\imath_\nu^*: \text{Gal} \left( \overline {\mathbb Q}_\nu / \mathbb Q_\nu \right) \hookrightarrow \text{Gal} \left( \overline {\mathbb Q} / \mathbb Q \right)$ given by $\sigma \mapsto \sigma \circ \imath_\nu$.  We consider $G_\nu = \text{Gal} \left( \overline {\mathbb Q}_\nu / \mathbb Q_\nu \right)$ a subgroup of the absolute Galois group $G_\mathbb Q = \text{Gal} \left( \overline {\mathbb Q} / \mathbb Q \right)$, and call it the decomposition group at $\nu$.  We define the inertia group $I_\nu$ as follows: When $\nu = p$ is a finite prime, $I_p$ is the kernel of the map $\text{Gal} \left( \overline {\mathbb Q}_p / \mathbb Q_p \right) \to \text{Gal} \left( \overline {\mathbb F}_p / \mathbb F_p \right)$ induced by reduction modulo $p$, and when $\nu = \infty$ is the infinite prime, $I_\infty = G_\infty = \text{Gal} \left( \mathbb C / \mathbb R \right)$ is the group generated by complex conjugation.  In either case, the inertia group is a normal subgroup of the decomposition group, and the quotient $G_\nu / I_\nu$ is a procyclic group with topological generator $\text{Frob}_\nu$.  For more properties, see \cite{MR82e:12016}.

Consider the restriction of $\rho$ to the inertia and decomposition groups.  We say that $\rho$ is unramified at $\nu$ if its restriction to the inertia group at $\nu$ is trivial, and ramified otherwise.  (We consider only $\ell$-adic representations which are ramified at finitely many places.)  We say $\rho$ is wildly ramified at a finite prime $\nu = p$ if its restriction to the Sylow pro-$p$-subgroup $I_p^w \subseteq I_p$ is nontrivial.  Note $\rho$ is ramified at the infinite prime $\nu = \infty$ if and only if $\det \rho(c) = -1$ on complex conjugation; we say that such representations are odd.  In this case,
\[ \rho(c) \simeq \begin{pmatrix} -1 \\ & 1 \end{pmatrix}. \]

\noindent When $\nu = \ell$ we say that $\rho$ is ordinary if the decomposition group may be conjugated over $\mathcal O$ to be upper-triangular with the upper entry on the diagional an unramified character i.e. there exist homomorphisms $\chi_0, \, \chi_\ell: G_\ell \to \mathcal O^\times$ such that
\[ \rho |_{G_\ell} \simeq \begin{pmatrix} \chi_\ell & * \\ & \chi_0 \end{pmatrix} \quad \text{with} \quad \chi_0 |_{I_\ell} = 1. \]

The cyclotomic character $\varepsilon_\ell: G_{\mathbb Q} \to {\mathbb Z}_\ell^\times$ is a canonical $\ell$-adic representation defined by acting on the $\ell$-power roots of unity: $\sigma \left( \zeta \right) = \zeta^{\varepsilon_\ell(\sigma)}$.  Note that $\varepsilon_\ell$ is unramified at $\nu \neq \ell, \, \infty$; in fact for $\ell$ odd, the Frobenius element $\text{Frob}_p \mapsto p$ for finite primes $p \neq \ell$ and complex conjugation $c \mapsto -1$ at the infinite prime, yet the inertia group $I_\ell$ has nontrivial image because there are no nontrivial $\ell$-power roots of unity in $\overline {\mathbb F}_\ell$.  In particular this gives a nontrivial action on $\ell$-power roots of unity $\mu_{\ell^\infty}$ so for any $k$-vector space $V$ acted upon by the absolute Galois group we define $V(1) = \mu_{\ell^\infty} \otimes V$ as the Tate twist of $V$.

\subsection{Adjoint Representation}

The residual representation $\overline \rho$ induces an action of the absolute Galois group on the $k$-vector space of $2 \times 2$ matrices given by $\sigma \cdot m = \overline \rho(\sigma) \ m \ \overline \rho(\sigma)^{-1}$; we denote this $k$-vector space with such an action by $\text{ad} \, \overline \rho$.  Explicitly, denote the matrices
\begin{equation} \label{adjoint_matrices} m_{-1} = \begin{pmatrix} 1 & 0 \\ 0 & 1 \end{pmatrix}, \quad m_0 = \begin{pmatrix} 0 & 0 \\ 1 & 0 \end{pmatrix}, \quad m_1 = \begin{pmatrix} 1 & 0 \\ 0 & -1 \end{pmatrix}, \quad m_2 = \begin{pmatrix} 0 & -1 \\ 0 & 0 \end{pmatrix}; \end{equation}

\noindent then $\text{ad} \, \overline \rho$ is the $k$-vector space spanned by the $m_j$.  We have a filtration
\[ \left \{ 0 \right \} = \text{ad}^3 \overline \rho \subset \text{ad}^2 \overline \rho \subset \text{ad}^1 \overline \rho \subset \text{ad}^0 \overline \rho \subset \text{ad}^{-1} \overline \rho = \text{ad} \, \overline \rho \]

\noindent where $\text{ad}^i \overline \rho$ is the $(3-i)$-dimensional subspace spanned by the $m_j$ for $j \geq i$.

\begin{infinitesimal} \label{infinitesimal} Let $\overline \rho: G_{\mathbb Q} \to GL_2(k)$ be a continuous Galois representation, and let $\epsilon$ be an infinitesimal i.e. $\epsilon^2 = 0$.  The equivalence classes of infinitesimal deformations $\overline \rho_\epsilon: G_{\mathbb Q} \to GL_2 \left( k[\epsilon] \right)$ satisfying $\overline \rho_\epsilon \equiv \overline \rho \mod {\epsilon \, k[\epsilon]}$ and $\det \overline \rho_\epsilon = \det \overline \rho$ are in one-to-one correspondence with the cohomology classes in $H^1 \left( G_\mathbb Q, \, \text{ad}^0 \overline \rho \right)$.
\end{infinitesimal}

\begin{proof}  We follow \cite[\S 4]{MR1638477}.  Given $\sigma \in G_{\mathbb Q}$ we may express an infinitesimal deformation in the form $\overline \rho_\epsilon(\sigma) = \left( 1_2 + \epsilon \, \xi_\sigma \right) \overline \rho(\sigma)$ for some $\xi_\sigma \in \text{Mat}_2(k)$, where $\xi_\sigma$ must have trace zero since $\det \overline \rho(\sigma) = \det \overline \rho_\epsilon(\sigma) = \left( 1 + \epsilon \, \text{tr} \, \xi_\sigma \right) \det \overline \rho(\sigma)$.  We claim that equivalence classes of homomorphisms are in one-to-one correspondence with $\xi \in H^1 \left( G_{\mathbb Q}, \, \text{ad}^0 \overline \rho \right)$.

For $\overline \rho_\epsilon$ to be a homomorphism we must have
\[ \begin{aligned} 1_2 & = \overline \rho_\epsilon (\sigma \, \tau) \cdot \overline \rho_\epsilon(\tau)^{-1} \cdot \overline \rho_\epsilon(\sigma)^{-1} \\ & = \left( 1_2 + \epsilon \, \xi_{\sigma \tau} \right) \overline \rho(\sigma \, \tau) \cdot \overline \rho(\tau)^{-1} \left( 1_2 - \epsilon \, \xi_\tau \right) \cdot \overline \rho(\sigma)^{-1} \left( 1_2 - \epsilon \, \xi_\sigma \right) \\ & = \left( 1_2 + \epsilon \, \xi_{\sigma \tau} \right) \cdot \left( 1_2 - \epsilon \left[ \overline \rho(\sigma) \ \xi_\tau \ \overline \rho(\sigma)^{-1} \right] \right) \cdot \left( 1_2 - \epsilon \, \xi_\sigma \right) \\ & = 1_2 + \epsilon \left[ \xi_{\sigma \tau} - \sigma \cdot \xi_\tau - \xi_\sigma \right]. \end{aligned} \]

\noindent  Hence $\overline \rho_\epsilon$ is a homomorphism if and only if $\xi$ is a cocycle i.e. $\xi_{\sigma \tau} = \xi_\sigma + \sigma \cdot \xi_\tau$.

The equivalence class of $\overline \rho_\epsilon$ consists of conjugates $A \, \overline \rho_\epsilon \, A^{-1}$ where $A = 1_2 + \epsilon \, \alpha$ for $\alpha \in \text{Mat}_2(k)$.  But $A \, \overline \rho_\epsilon(\sigma) \, A^{-1} = \left( 1_2 + \epsilon \left[ \alpha - \sigma \cdot \alpha \right] \right) \overline \rho_\epsilon(\sigma)$ so equivalence classes yield $\xi_\sigma$ modulo coboundaries i.e. $\alpha - \sigma \cdot \alpha$.  Hence $\xi \in H^1 \left( G_{\mathbb Q}, \, \text{ad}^0 \overline \rho \right)$ as claimed.  \end{proof}

\subsection{Selmer Groups}

We wish to classify infinitesimal deformations even further by considering what happens when we restrict to the decomposition and inertia groups.  Recall that for each place $\nu$ of $\mathbb Q$, we have restriction maps $\text{res}_\nu: H^1 \left( G_{\mathbb Q}, \, \text{ad}^0 \overline \rho \right) \to H^1 \left( G_\nu, \, \text{ad}^0 \overline \rho \right)$ as well as the ``inflation-restriction'' exact sequence
\[ 0 \to H^1 \left( G_\nu / I_\nu, \, \left( \text{ad}^0 \overline \rho \right)^{I_\nu} \right) \to H^1 \left( G_\nu, \, \text{ad}^0 \overline \rho \right) \to H^1 \left( I_\nu, \, \text{ad}^0 \overline \rho \right). \]

\noindent  For more properties, see \cite{MR88e:14028}.

We will use these maps to define subgroups $H^1_f \left( G_\nu, \, \text{ad}^0 \overline \rho \right) \subseteq H^1 \left( G_\nu, \, \text{ad}^0 \overline \rho \right)$ that will encode information about infinitesimal deformations.  When $\nu \neq \ell$ each deformation $\overline \rho_\epsilon$ should be (un)ramified when $\overline \rho$ is (un)ramified so the restriction of a class from $H^1 \left( G_\nu, \, \text{ad}^0 \overline \rho \right)$ to $H^1 \left( I_\nu, \, \text{ad}^0 \overline \rho \right)$ should be trivial.  Define
\[ \begin{aligned} H^1_f \left( G_\nu, \, \text{ad}^0 \overline \rho \right) & = \text{ker} \left[H^1 \left( G_\nu, \, \text{ad}^0 \overline \rho \right) \to H^1 \left( I_\nu, \, \text{ad}^0 \overline \rho \right) \right] \\ & = H^1 \left( G_\nu / I_\nu, \, \left( \text{ad}^0 \overline \rho \right)^{I_\nu} \right)  \end{aligned} \qquad \text{for $\nu \neq \ell$.} \]

\noindent When $\nu = \ell$ and $\overline \rho$ is ordinary i.e. the restriction of $\overline \rho$ to $G_\ell$ is upper-triangular, each deformation $\overline \rho_\epsilon$ should be ordinary as well.  We choose
\[ \begin{aligned} H^1_f \left( G_\ell, \, \text{ad}^0 \overline \rho \right) & \subseteq \text{ker} \left[ H^1 \left( G_\ell, \, \text{ad}^0 \overline \rho \right) \to H^1 \left( G_\ell, \, \text{ad}^0 \overline \rho / \text{ad}^1 \overline \rho \right) \right] \\ & = \text{im} \left[ H^1 \left( G_\ell, \, \text{ad}^1 \overline \rho \right) \to H^1 \left( G_\ell, \, \text{ad}^0 \overline \rho \right) \right]. \end{aligned} \]
\noindent When restricted to the inertia group, the diagonal terms of $\overline \rho_\epsilon$ should be the same as those of $\overline \rho$ so the restriction of a class from $H^1 \left( G_\ell, \, \text{ad}^1 \overline \rho \right)$ to $H^1 \left( I_\ell, \, \text{ad}^1 \overline \rho / \text{ad}^2 \overline \rho \right)$ should be trivial.  Define
\begin{equation} \label{local_condition} \begin{aligned} H^1_f & \left( G_\ell, \, \text{ad}^0 \overline \rho \right) \\ & = \text{im} \biggl[ \text{ker} \left[ H^1 \left( G_\ell, \, \text{ad}^1 \overline \rho \right) \to H^1 \left( I_\ell, \, \text{ad}^1 \overline \rho / \text{ad}^2 \overline \rho \right) \right] \to H^1 \left( G_\ell, \, \text{ad}^0 \overline \rho \right) \biggr]. \end{aligned} \end{equation}

\noindent We explain how this definition for $\nu = \ell$ compares to the usual one.  In general we have the exact sequence
\[ \begin{CD} H^1_f \left( G_\ell, \, \text{ad}^0 \overline \rho \right) @>>> H^1 \left( G_\ell, \, \text{ad}^0 \overline \rho \right) @>>> H^1 \left( I_\ell, \, \text{ad}^0 \overline \rho / \text{ad}^2 \overline \rho \right)  \end{CD} \]

\noindent but the map $H^1 \left( G_\ell, \, \text{ad}^1 \overline \rho \right) \to H^1 \left( G_\ell, \, \text{ad}^0 \overline \rho \right)$ is not an injection.  However when $\det \rho= \varepsilon_\ell$ is the cyclotomic character, the group $H^0 \left( G_\ell, \,  \text{ad}^0 \overline \rho / \text{ad}^1 \overline \rho \right)$ is trivial, so we recover the usual definition as considered in \cite{MR96d:11071} and \cite{MR99d:11067b}:
\[ H^1_f \left( G_\ell, \, \text{ad}^0 \overline \rho \right) = \text{ker} \left[ H^1 \left( G_\ell, \, \text{ad}^0 \overline \rho \right) \to H^1 \left( I_\ell, \, \text{ad}^0 \overline \rho / \text{ad}^2 \overline \rho \right) \right]. \]

Fix a finite set $\Sigma$ of places that does not contain $\ell$.  We define the Selmer group for $\text{ad}^0 \overline \rho$ with respect to $\Sigma$ as
\[ H^1_\Sigma \left( \mathbb Q, \, \text{ad}^0 \overline \rho \right) = \text{ker} \left[ H^1 \left( G_{\mathbb Q}, \, \text{ad}^0 \overline \rho \right) \to \bigoplus_{\nu \not \in \Sigma} H^1 \left( G_\nu, \, \text{ad}^0 \overline \rho \right) \left/ H^1_f \left( G_\nu, \, \text{ad}^0 \overline \rho \right) \right. \right] \]

\noindent i.e. the collection of classes $\xi \in H^1 \left( G_{\mathbb Q}, \, \text{ad}^0 \overline \rho \right)$ such that $\text{res}_\nu(\xi) \in H^1_f \left( G_\nu, \, \text{ad}^0 \overline \rho \right)$ for all places $\nu \not \in \Sigma$.

Define $H^1_f \left( G_\nu, \, \text{ad}^0 \overline \rho(1) \right)$ as the annihilator of $ H^1_f \left( G_\nu, \, \text{ad}^0 \overline \rho \right)$ under the pairing
\[ H^1 \left( G_\nu, \, \text{ad}^0 \overline \rho \right) \times H^1 \left( G_\nu, \, \text{ad}^0 \overline \rho(1) \right) \to H^2 \left( G_\nu, \, k(1) \right) \simeq k \]

\noindent and define the dual Selmer group $H^1_\Sigma \left( \mathbb Q, \, \text{ad}^0 \overline \rho(1) \right)$ as classes $\xi \in H^1 \left( G_{\mathbb Q}, \, \text{ad}^0 \overline \rho(1) \right)$ such that $\text{res}_\nu(\xi) \in H^1_f \left( G_\nu, \, \text{ad}^0 \overline \rho(1) \right)$ for all places $\nu \not \in \Sigma$, yet $\text{res}_\nu(\xi) = 0$ otherwise.

\begin{goins_lemma} \label{goins_lemma} For $\ell$ an odd prime, let $\overline \rho: G_\ell \to GL_2(k)$ be a continuous, ordinary  mod $\ell$ representation that is wildly ramified. Then \[ \dim_k H^1_f \left( G_\ell, \, \text{ad}^0 \overline \rho \right) = 1 + \dim_k H^0 \left( G_\ell, \, \text{ad}^0 \overline \rho \right). \]
\end{goins_lemma}

\begin{proof} As $\overline \rho$ is ordinary, assume that a given $\sigma \in G_\ell$ maps to
\[ \overline \rho(\sigma) = \begin{pmatrix} a & b \\ & d \end{pmatrix} \quad \text{for some} \quad a, \, d \in k^\times, \quad b \in k. \]

\noindent Recall that $\text{ad}^i \overline \rho$ is spanned by $m_j$ as in \eqref{adjoint_matrices} for $j \geq i$, and note the identities
\[ \sigma \cdot m_0 = \frac da \, m_0 + \frac ba \, m_1 + \frac {b^2}{a \, d} \, m_2, \quad \sigma \cdot m_1 = m_1 + \frac {2 \, b}d \, m_2, \quad \sigma \cdot m_2 = \frac ad \, m_2. \]

\noindent In particular, $G_\ell$ acts trivially on $V = \text{ad}^1 \overline \rho / \text{ad}^2 \overline \rho$.  To prove the proposition, note
\[ H^1_f \left( G_\ell, \, \text{ad}^0 \overline \rho \right) \simeq \frac {K_1}{K_1 \cap K_2} \quad \text{where} \quad \begin{aligned} K_1 & = \text{ker} \left[H^1 \left( G_\ell, \, \text{ad}^1 \overline \rho \right) \to H^1 \left( I_\ell, \, V \right) \right], \\ K_2 & = \text{ker} \left[ H^1 \left( G_\ell, \, \text{ad}^1 \overline \rho \right) \to H^1 \left( G_\ell, \, \text{ad}^0 \overline \rho \right) \right]; \end{aligned} \]

\noindent we will compute the dimension of $K_1$ and show $K_1 \cap K_2$ is trivial.

As $\overline \rho$ is wildly ramified, there exists $\sigma_0 \in I_\ell^w$ (wild inertia) such that
\[ \overline \rho(\sigma_0) = \begin{pmatrix} 1 & b_0 \\ & 1 \end{pmatrix} \quad \text{for some  $b_0 \in k^\times$.} \]

\noindent When $j > i$ the elements of $\text{ad}^i \overline \rho / \text{ad}^j \overline \rho$ fixed by $\sigma_0$ form a 1-dimensional space that only depends on $j$: 
\[ H^0 \left( G_\ell, \, \text{ad}^i \overline \rho / \text{ad}^j \overline \rho \right) \subseteq H^0 \left( I_\ell, \, \text{ad}^i \overline \rho / \text{ad}^j \overline \rho \right) \subseteq H^0 \left( I_\ell^w, \, \text{ad}^i \overline \rho / \text{ad}^j \overline \rho \right) = \text{ad}^{j-1} \overline \rho / \text{ad}^j \overline \rho. \]

We compute the dimension of $K_1$, the kernel of the composition $H^1 \left( G_\ell, \text{ad}^1 \overline \rho \right) \to H^1 \left( G_\ell, V \right) \to H^1 \left( I_\ell, V \right)$.  Wild inertia acts trivially on $\text{ad}^2 \overline \rho$ and hence nontrivially on $\text{ad}^2 \overline \rho(1)$ so $H^2 \left( G_\ell, \, \text{ad}^2 \overline \rho \right) \simeq H^0 \left( G_\ell, \, \text{ad}^2 \overline \rho(1) \right) \subseteq H^0 \left( I_\ell^w, \, \text{ad}^2 \overline \rho(1) \right)$ is trivial. Then the exact sequence $0 \to \text{ad}^2 \overline \rho \to \text{ad}^1 \overline \rho \to V \to 0$ implies the exact sequence
\[ 0 \to H^0 \left( G_\ell, \, V \right) \to H^1 \left( G_\ell, \, \text{ad}^2 \overline \rho \right) \to H^1 \left( G_\ell, \, \text{ad}^1 \overline \rho \right) \to H^1 \left( G_\ell, \, V \right) \to 0. \]

\noindent In particular, the map $H^1 \left( G_\ell, \text{ad}^1 \overline \rho \right) \to H^1 \left( G_\ell, V \right)$ is surjective.  On the other hand, the map $H^1 \left( G_\ell, V \right) \to H^1 \left( I_\ell, V \right)$ has kernel $H^1 (G_\ell/I_\ell, \, V^{I_\ell}) \simeq \text{Hom}(k, \, k)$ of dimension 1 because $G_\ell/I_\ell$ is a procyclic group with trivial action on $V$.  Hence
\[ \begin{aligned} \dim_k K_1 & = \dim_k \ker \left[ H^1 \left( G_\ell, \text{ad}^1 \overline \rho \right) \to H^1 \left( G_\ell, V \right) \right] \\ & \qquad + \dim_k \ker \left[ H^1 \left( G_\ell, V \right) \to H^1 \left( I_\ell, V \right) \right] \\ & = \left[ \dim_k H^1 \left( G_\ell, \, \text{ad}^2 \overline \rho \right) - \dim_k H^0 \left( G_\ell, \, V \right) \right] + \dim_k H^1 \left( G_\ell/I_\ell, \, V^{I_\ell} \right) \\ & = \dim_k H^1 \left( G_\ell, \, \text{ad}^2 \overline \rho \right) = 1 + \dim_k H^0 \left( G_\ell, \, \text{ad}^2 \overline \rho \right) \\ & = 1 + \dim_k H^0 \left( G_\ell, \, \text{ad}^0 \overline \rho \right) . \end{aligned} \]

We show $K_1 \cap K_2$ is trivial.  As $H^0 (G_\ell, \text{ad}^1 \overline \rho) = H^0 (G_\ell, \text{ad}^0 \overline \rho)$ we have the exact sequence $0 \to H^0 \left( G_\ell, \, \text{ad}^0 \overline \rho / \text{ad}^1 \overline \rho \right) \to H^1 \left( G_\ell, \, \text{ad}^1 \overline \rho \right) \to H^1 \left( G_\ell, \, \text{ad}^0 \overline \rho \right)$, so
\[ K_2 \simeq H^0 \left( G_\ell, \, \text{ad}^0 \overline \rho / \text{ad}^1 \overline \rho \right) \subseteq H^0 \left( I_\ell, \, \text{ad}^0 \overline \rho / \text{ad}^1 \overline \rho \right). \]

\noindent Either $K_2$ is 1-dimensional and we have equality or $K_2$ is trivial.  A brief chase of the exact diagram
\[ \begin{CD} @. 0 @. 0 @. 0 \\ @. @VVV @VVV @VVV \\ 0 @>>> K_1 \cap K_2 @>>> K_2 @>>> H^0 \left( I_\ell, \, \text{ad}^0 \overline \rho / \text{ad}^1 \overline \rho \right) \\ @. @VVV @VVV @VVV \\ 0 @>>> K_1 @>>> H^1 \left( G_\ell, \, \text{ad}^1 \overline \rho \right) @>>> H^1 \left( I_\ell, \, \text{ad}^1 \overline \rho / \text{ad}^2 \overline \rho \right) \\ @. @VVV @VVV @VVV \\ @. H^1_f \left( G_\ell, \, \text{ad}^0 \overline \rho \right) @>>> H^1 \left( G_\ell, \, \text{ad}^0 \overline \rho \right) @>>> H^1 \left( I_\ell, \, \text{ad}^0 \overline \rho / \text{ad}^2 \overline \rho \right) \end{CD} \]

\noindent  shows $K_1 \cap K_2$ is trivial in either case. \end{proof}

\subsection{Deformations of type $\Sigma$}

Again, fix a finite set $\Sigma$ of places that does not contain $\ell$.  Assume that $\rho: G_{\mathbb Q} \to GL_2(\mathcal O)$ is ordinary and ramified at finitely many primes, and its residual representation $\overline \rho: G_{\mathbb Q} \to GL_2(k)$ is absolutely irreducible and wildly ramified at $\ell$.  Denote by $\mathcal O'$ a complete, local Noetherian $\mathcal O$-algebra with field of fractions $K'$ having residue field $k' = \mathcal O' / \lambda'$.  We say a representation $\rho': G_{\mathbb Q} \to GL_2(\mathcal O')$ is a deformation of $\overline \rho$ of type $\Sigma$ if
\begin{enumerate}
\item $\overline \rho' \simeq \overline \rho \otimes_k k'$ and $\det \rho' = \det \rho$;
\item $\rho' |_{I_\nu} \simeq \rho |_{I_\nu} \otimes_K K'$ for all $\nu \not \in \Sigma$; and
\item $\rho' |_{G_\ell} \simeq \begin{pmatrix} \chi_\ell & * \\ & \chi_0 \end{pmatrix}$ where $\chi_\ell = \det \rho \cdot {\chi_0}^{-1}$ and $\chi_0 |_{I_\ell} = 1$.
\end{enumerate}

\noindent These conditions reflect via Proposition \ref{infinitesimal} precisely the cohomology classes in the Selmer group $H_\Sigma^1 \left( \mathbb Q, \, \text{ad}^0 \overline \rho \right)$.  Indeed, an equivalence class of deformations $\rho'$ satisfying condition (1) corresponds to $\xi \in H^1 \left( G_\mathbb Q, \, \text{ad}^0 \overline \rho \right)$; condition (2) corresponds to $\text{res}_\nu(\xi) \in H_f^1 \left( G_\nu, \, \text{ad}^0 \overline \rho \right)$; and condition (3) corresponds to $\text{res}_\ell(\xi) \in H_f^1 \left( G_\ell, \, \text{ad}^0 \overline \rho \right)$.

The local conditions above are also reflected in the Artin conductor.  Denote the conductor of $\rho$ by $N_\emptyset$.  The condition at the infinite place $\nu = \infty$ is equivalent to saying a deformation $\rho'$ is odd precisely when $\rho$ is odd, so consider the finite places $\nu = p$.  If $\rho$ is unramified at $p \not \in \Sigma$ then so is $\rho'$; in particular $\rho'$ is ramified at finitely many primes.  For $p \in \Sigma$ there are no conditions on the restriction of $\rho'$ at $I_p$ so at worst the $p$-part of the conductor is $p^2$.  Hence the conductor of $\rho'$ is divisible by $N_\emptyset$ yet divides $N_\Sigma = N_\emptyset \, \prod_p p^2$ in terms of the product over $p \in \Sigma$ not dividing $N_\emptyset$.  In particular, when = $N_\emptyset = N(\overline \rho) \, \ell$ in terms of the Serre conductor defined in \cite{MR93h:11124}, the local conditions above force a deformation $\rho'$ to be minimally ramified in the sense of \cite[Definition 3.1]{MR1638490}.

The following proposition shows that the ``deformation conditions'' above are ``representable'' i.e. there exists a universal deformation ring $R_\Sigma$. 

\begin{existence} \label{existence} Let $\rho: G_{\mathbb Q} \to GL_2(\mathcal O)$ be a continuous $\ell$-adic representation such that $\rho$ is ordinary and ramified at finitely many primes, while $\overline \rho$ is absolutely irreducible and wildly ramified at $\ell$.  Fix a finite set $\Sigma$ of places that does not contain $\ell$.  There exists a universal deformation $\rho_\Sigma^{\text{univ}}: G_{\mathbb Q} \to GL_2(R_\Sigma)$ of $\overline \rho$ of type $\Sigma$ in the sense that if $\rho': G_{\mathbb Q} \to GL_2(\mathcal O')$ is any deformation of $\overline \rho$ of type $\Sigma$, then there exists a unique $\mathcal O$-algebra homomorphism $\phi: R_\Sigma \to \mathcal O'$ such that $\rho' \simeq \phi \circ \rho_\Sigma^{\text{univ}}$.
\end{existence}

It may appear that such a result follows from the well-known properties of deformation rings associated to absolutely irreducible residual representations, but we do not impose the typical local condition that $\chi_\ell$ is ramified at $\ell$, as in \cite[p. 457]{MR96d:11071} and \cite[\S 30]{MR1638481}.  Indeed, if $\chi_\ell$ and $\chi_0$ were trivial when restricted to the decomposition group at $\ell$ -- and hence both unramified at $\ell$ -- then the centralizer of the residual representation $\overline \rho \vert_{G_\ell}$ would be non-scalar, and hence a desired local ``universal'' deformation $G_\ell \to GL_2(R_\Sigma)$ constructed out of the local representation would be ``versal'' at best.  We exploit the fact that the residual representation is wildly ramified at $\ell$ to produce a global universal deformation $G_{\mathbb Q} \to GL_2(R_\Sigma)$.

\begin{proof} We follow the exposition in \cite{MR1638481}.  Consider the functor
\[ \mathcal D_\Sigma: \quad \left \{ \begin{matrix} \text{artinian $\mathcal O$-algebras $A$} \\ \text{having residue field $k$} \\ \text{with a local $\mathcal O$-homomorphism $A \to \mathcal O$} \end{matrix} \right \} \to \left \{ \begin{matrix} \text{equivalence classes} \\ \text{of deformations $\rho_A$} \\ \text{of $\overline \rho$ of type $\Sigma$} \end{matrix} \right \}. \]

\noindent We say this functor is representable if there exists a complete, local Noetherian $\mathcal O$-algebra $R_\Sigma$ depending only on the representation $\rho$ and the set $\Sigma$ such that
\[ \mathcal D_\Sigma \left( A \right) \simeq \text{Hom}_{\mathcal O} \left( R_\Sigma, \, A \right). \]

\noindent If this is the case, then the image $\mathcal D_\Sigma \left( R_\Sigma \right)$ associated with the trivial homomorphism $R_\Sigma \to R_\Sigma$ corresponds to a deformation $\rho_\Sigma^{\text{univ}}$ of $\overline \rho$ of type $\Sigma$.  Moreover, the isomorphism above implies that any image $\mathcal D_\Sigma (A)$ may be associated with a unique $\mathcal O$-algebra homomorphism $\phi_A: R_\Sigma \to A$, so that $\rho_A \simeq \phi_A \circ \rho_\Sigma^{\text{univ}}$.  When $A = \mathcal O' / \Lambda'^n$ we have the inverse limit
\[ \mathcal O' = \projlim_{n \to \infty} \mathcal O' / \Lambda'^n \implies \rho' \simeq \phi \circ \rho_\Sigma^{\text{univ}} \quad \text{where} \quad \phi: R_\Sigma \to \mathcal O', \]
\noindent so that the proposition follows.  Hence it suffices to show that $\mathcal D_\Sigma$ is representable.

Now consider the functor $\mathcal D_{\overline \rho}$ sending artinian $\mathcal O$-algebras $A$ to deformations $\rho_A$ satisfying just the deformation condition (1) above.  It is well-known that this functor is representable because $\overline \rho$ is absolutely irreducible and ramified at finitely many primes; see \cite[\S1.2, Proposition 1]{MR90k:11057} and \cite[\S 1.3, (a.3)]{MR90k:11057}.  To show $\mathcal D_\Sigma$ is representable, we use \cite[\S 26, Proposition 1]{MR1638481} and \cite[\S 23, Corollary]{MR1638481} to note that it suffuces to show that the local conditions in (2) and (3) are ``deformation conditions'' in the sense of \cite[\S 23]{MR1638481} defined as follows:  Fix a place $\nu$ of $\mathbb Q$ and a free rank two $\mathcal O$-module $V$ such that $\rho: G_\nu \to GL(V)$.  Consider a category whose objects are pairs $(A,V_A)$ of artinian $\mathcal O$-algebras $A$ along with free rank two $A$-modules $V_A$ endowed with a continuous action by $G_\nu$.  Note that we may consider an object to be an assignment $A \rightsquigarrow \rho_A$ sending an artinian $\mathcal O$-algebra $A$ to a 2-dimensional representaton $\rho_A: G_\nu \to GL(V_A)$ defined by this action.  A morphism $(A,V_A) \to (C, V_C)$ in this category consists of an $\mathcal O$-algebra homomorphism $A \to C$ along with an $A$-module homomorphism $V_A \to V_C$ inducing the isomorphism $V_C \simeq V_A \otimes_A C$.  We say that a category $\mathcal D_\nu$ is a deformation condition for $\overline \rho$ if the following four properties hold for all objects:
\begin{enumerate}
\item[DC1:] $(k, V \otimes_{\mathcal O} k) \in \mathcal D_\nu$ where the action on $V \otimes_{\mathcal O} k$ by $G_\nu$ is given by $\overline \rho \vert_{G_\nu}$.
\item[DC2:] For any morphism $(A,V_A) \to (C,V_C)$, if $(A,V_A) \in \mathcal D_\nu$ then $(C,V_C) \in \mathcal D_\nu$.
\item[DC3:] For any morphism $(A,V_A) \to (C,V_C)$, if $(C,V_C) \in \mathcal D_\nu$ and $A \hookrightarrow C$ then $(A,V_A) \in \mathcal D_\nu$.
\item[DC4:] Say that the following diagram of morphisms is cartesian:
\[ \begin{matrix} D \\  \swarrow \quad \searrow \\ A \qquad \qquad B \\ \alpha \quad \searrow \quad \swarrow \quad \beta \\ C \end{matrix} \]
\noindent that is, $D \simeq A \times_C B$ is the fiber product.  The object $(D,V_D) \in \mathcal D_\nu$ if and only if both $(A,V_A), (B,V_B) \in \mathcal D_\nu$.
\end{enumerate}

For each place $\nu$, consider the full subcategory $\mathcal D_\nu \subseteq \mathcal D_{\overline \rho}$ of objects $(A,V_A)$ such that the action by inertia is given by $\rho \vert_{I_\nu} \otimes_{\mathcal O} A$.  It is clear that a deformation $\rho_A$ of $\overline \rho$ satisfies conditions (1) and (2) if and only if $(A,V_A)$ is in $\mathcal D_\nu$ for all $\nu \not \in \Sigma$.  We claim that a deformation $\rho_A$ of $\overline \rho$ satisfies conditions (1) and (3) if and only if $(A,V_A)$ is in $\mathcal D_\ell$.  To this end, choose an object $(A,V_A) \in \mathcal D_\ell$, and consider the $A$-linear submodule
\[ (V_A)_{I_\ell} = \biggl \{ v \in V_A \, \biggl \vert \, \text{$\sigma \cdot v = \det \rho(\sigma) \cdot v$ for all $\sigma \in I_\ell$} \biggr. \biggr \}. \]

\noindent This submodule is nontrivial because there exists an $A$-basis of $V_A$ such that
\[ \rho_A \vert_{I_\ell} \simeq \rho \vert_{I_\ell} \otimes_{\mathcal O} A \simeq \begin{pmatrix} \det \rho & * \\ & 1 \end{pmatrix}; \]

\noindent yet the submodule is not all of $V_A$ because $\rho$ is wildly ramified at $\ell$ i.e. $\ast$ is not identically zero.  We identify $(V_A)_{I_\ell} \simeq V_{I_\ell} \otimes_{\mathcal O} A$ as the free rank one $A$-module of $I_\ell$-coinvariants because inertia acts trivially on the quotient $(V_A)^{I_\ell} = V_A / (V_A)_{I_\ell} \simeq V^{I_\ell} \otimes_{\mathcal O} A$ of $I_\ell$-invariants.  We mention in passing that
\[ \text{ad} \, \overline \rho \simeq \text{Hom}_k \left( V \otimes_{\mathcal O} k, \, V \otimes_{\mathcal O} k \right) \quad \text{and} \quad \text{ad}^2 \, \overline \rho \simeq \text{Hom}_k \left( V_{I_\ell} \otimes_{\mathcal O} k, \, V^{I_\ell} \otimes_{\mathcal O} k \right). \]

\noindent Since $I_\ell$ is a normal subgroup of $G_\ell$, the $A$-submodules $(V_A)_{I_\ell}$ and $(V_A)^{I_\ell}$ are stable under action by $G_\ell$.   To see why, choose $\tau \in G_\ell$ and $v \in (V_A)_{I_\ell}$; we must show that $\sigma \cdot (\tau \cdot v) = \det \rho(\sigma) \cdot (\tau \cdot v)$ for all $\sigma \in I_\ell$.  Indeed, we have
\[ \sigma \cdot \tau \cdot v = \tau \cdot \sigma' \cdot v = \tau \cdot \det \rho(\sigma') \cdot v = \tau \cdot \det \rho(\sigma) \cdot v = \det \rho(\sigma) \cdot \tau \cdot v \]

\noindent where $\sigma' = \tau^{-1} \cdot \sigma \cdot \tau \in I_\ell$.  This shows that $\rho_A \vert_{G_\ell}$ acts on $(V_A)_{I_\ell}$ and $(V_A)^{I_\ell}$ by some characters $\chi_\ell$ and $\chi_0$, respectively.  It is clear that $\chi_0 \cdot \chi_\ell = \det \rho_A = \det \rho$ and $\chi_0$ is unramified.  Hence $\rho_A$ is a deformation of $\overline \rho$ of type $\Sigma$ is and only if $(A,V_A)$ is in $\mathcal D_\nu$ for all places either $\nu \not \in \Sigma$ or $\nu = \ell$.

It suffices show that $\mathcal D_\nu$ is a deformation condition for $\overline \rho$ for all places $\nu$.  It is clear that $(k, V \otimes_{\mathcal O} k)$ is in $\mathcal D_\nu$ i.e. (DC1) is satisfied.  Say that we have a morphism $(A,V_A) \to (C,V_C)$ where $(A,V_A)$ is in $\mathcal D_\nu$.  Since $V_C \simeq V_A \otimes_A C$ as $A$-modules, the action by inertia must be given by $\left( \rho \vert_{I_\nu} \otimes_{\mathcal O} A \right) \otimes_A C \simeq \rho \vert_{I_\nu} \otimes_{\mathcal O} C$; hence $(C, V_C)$ is in $\mathcal D_\nu$ as well.  This shows that (DC2) is satisfied.  Say that we have a morphism $(A,V_A) \to (C,V_C)$ where $(C,V_C)$ is in $\mathcal D_\nu$ and we have the exact sequence of $A$-modules:
\[ \begin{CD} 0 @>>> A @>>> C @>>> C/A @>>> 0. \end{CD} \]

\noindent We tensor this sequence over $A$ with the free -- and hence flat -- $A$-module $V_A$ to find the exact sequence
\[ \begin{CD} 0 @>>> V_A @>>> V_C @>>> V_A \otimes_A C/A @>>> 0. \end{CD} \]

\noindent Fix $\sigma \in I_\nu$ and $v \in V_A$ and let $w = \left( \sigma - \rho(\sigma) \right) \cdot v \in V_A$.  By assumption on the object $(C,V_C)$ we have $w \in \text{ker} \left( V_A \to V_C \right)$. However $V_A$ injects into $V_C$ so $w = 0$ i.e. the action of inertia on $V_A$ is given by $\rho \vert_{I_\nu} \otimes_{\mathcal O} A$.  This shows that (DC3) is satisfied.

To show (DC4) is satisfied, assume that we have $\mathcal O$-algebra homomorphisms $\alpha: A \to C$ and $\beta: B \to C$, and denote $D = A \times_C B$ as the fiber product.  Fix a free 2-dimensional $D$-module $V_D$, and denote $V_A$ and $V_B$ as the tensor products of $V_D$ with respect to the projections $1 \times_C \beta$ and $\alpha \times_C 1$, respectively; then $V_D \simeq V_A \times_{V_C} V_B$.  If $(D,V_D)$ is in $\mathcal D_\nu$, then the canonical morphisms $(D,V_D) \to (A,V_A)$ and $(D,V_D) \to (B,V_B)$ show that both $(A,V_A)$ and $(B,V_B)$ are in $\mathcal D_\nu$.  Conversely, if the action of inertia on both $V_A$ and $V_B$ are specified, then it is also specified on $V_C \simeq V_A \otimes_A C \simeq V_B \otimes_B C$ and hence it is specified on $V_D \simeq V_A \times_{V_C} V_B$ as well.  This shows that if $(A,V_A)$ and $(B,V_B)$ are in $\mathcal D_\nu$ then so is $(D,V_D)$.  \end{proof}

\subsection{Dimension Computations}

We conclude this section by computing the dimension of the universal deformation ring.  Note that this is the only proposition in this section which assumes the residual representation is odd.

\begin{dimension} \label{dimension} For $\ell$ an odd prime, let $\rho: G_{\mathbb Q} \to GL_2(\mathcal O)$ be a continuous $\ell$-adic representation such that $\rho$ is ordinary and ramified at finitely many primes, while $\overline \rho$ is odd, absolutely irreducible, and wildly ramified at $\ell$.  If $\ell = 3$ assume moreover that $\overline \rho$ remains absolutely irreducible when restricted to $\text{Gal} \left( \overline {\mathbb Q} / \mathbb Q(\sqrt{-3}) \right)$.
\begin{enumerate}
\item Fix a finite set $\Sigma$ of places that does not contain $\ell$.  The universal deformation ring $R_\Sigma$ of $\overline \rho$ can be topologically generated as an $\mathcal O$-algebra by $\dim_k H^1_\Sigma \left( \mathbb Q, \, \text{ad}^0 \overline \rho \right)$ elements.
\item  Fix a finite set $\Sigma$ of finite places $q \equiv 1 \mod \ell$ such that $\overline \rho$ is unramified at $q$ and $\overline \rho(\text{Frob}_q)$ has distinct $k$-rational eigenvalues.  Then
\[ \dim_k H^1_\Sigma \left( \mathbb Q, \, \text{ad}^0 \overline \rho \right) = \# \Sigma + \dim_k H^1_\Sigma \left( \mathbb Q, \, \text{ad}^0 \overline \rho(1) \right). \]
\end{enumerate}
\end{dimension}

\begin{proof}  The Selmer group of $\text{ad}^0 \overline \rho$ with respect to $\Sigma$ may also be thought of as the kernel of the map from $H^1 \left( G_{\mathbb Q}, \, \text{ad}^0 \overline \rho \right)$ to $\bigoplus_\nu H^1 \left( G_\nu, \, \text{ad}^0 \overline \rho \right) / L_\nu$ where $L_\nu = H^1_f \left( G_\nu, \, \text{ad}^0 \overline \rho \right)$ when $\nu \not \in \Sigma$ and $L_q = H^1 \left( G_q, \, \text{ad}^0 \overline \rho \right)$ otherwise.  In order to prove the first statement, it suffices to exhibit an isomorphism
\begin{equation} \label{mazur} \text{Hom}_k \left( \mathfrak m_\Sigma / \left( \lambda, \, { \mathfrak m_\Sigma}^2 \right), \, k \right) \simeq H^1_\Sigma \left( \mathbb Q, \, \text{ad}^0 \overline \rho \right) \end{equation}

\noindent  with $\mathfrak m_\Sigma$ the maximal ideal of $R_\Sigma$.  A $k$-linear homomorphism $\mathfrak m_\Sigma / \left( \lambda, \, { \mathfrak m_\Sigma}^2 \right) \to k$ defines an $\mathcal O$-algebra homomorphism $R_\Sigma \to k[\epsilon]$ (and vice-versa) so we find a deformation $G_{\mathbb Q} \to GL_2(R_\Sigma) \to GL_2(k[\epsilon])$ which gives a class in $\xi \in H^1 \left( G_{\mathbb Q}, \, \text{ad}^0 \overline \rho \right)$.  It is clear from the construction that $\text{res}_\nu(\xi) \in L_\nu$ for all places $\nu$ so we have the desired map.

We now prove the second statement.  By \cite{MR96d:11071} we have the identity
\begin{equation} \label{wiles_ratio} \frac {\# H^1_\Sigma \left( \mathbb Q, \, \text{ad}^0 \overline \rho \right)}{\# H^1_\Sigma \left( \mathbb Q, \, \text{ad}^0 \overline \rho(1) \right)} = \frac {\# H^0 \left( G_{\mathbb Q}, \, \text{ad}^0 \overline \rho \right)}{\# H^0 \left( G_{\mathbb Q}, \, \text{ad}^0 \overline \rho(1) \right)} \prod_\nu \frac {\# L_\nu}{\# H^0 \left( G_\nu, \, \text{ad}^0 \overline \rho \right)}. \end{equation}

First consider $\nu \not \in \Sigma$.  When $\nu = p \neq \ell$ is a finite prime we have $\# L_p = \# H^1 \left( G_p / I_p, \, \left( \text{ad}^0 \overline \rho \right)^{I_p} \right) = \# H^0 \left( G_p, \, \text{ad}^0 \overline \rho \right)$.  When $\nu = \infty$ is the infinite prime $\# L_\infty = 1$ since $G_\infty = I_\infty$ so consider $H^0 \left( \text{Gal}\left( \mathbb C / \mathbb R \right), \, \text{ad}^0 \overline \rho \right)$ i.e. the elements fixed by complex conjugation.  By assumption $\overline \rho$ is odd so $c \cdot m_j = (-1)^{j-1} m_j$ with notation as in \eqref{adjoint_matrices}; hence $\# H^0 \left( G_\infty, \, \text{ad}^0 \overline \rho \right) = \# k$.

Now consider $q \in \Sigma$.  By the Local Duality theorems,
\[ \frac {\# H^1 \left( G_q, \, \text{ad}^0 \overline \rho \right)}{\# H^0 \left( G_q, \, \text{ad}^0 \overline \rho \right)} = \# H^2 \left( G_q, \, \text{ad}^0 \overline \rho \right) = \# H^0 \left( G_q, \, \text{ad}^0 \overline \rho(1) \right). \]

\noindent As $q \equiv 1 \pod \ell$ the decomposition group $G_q$ acts trivially on the $\ell$-power roots of unity: $H^0 \left( G_q, \, \text{ad}^0 \overline \rho(1) \right) = H^0 \left( G_q, \, \text{ad}^0 \overline \rho \right)$.  Explicitly denote
\[ \overline \rho \left( \text{Frob}_q \right) = \begin{pmatrix} \alpha_q \\ & \beta_q \end{pmatrix} \implies \text{Frob}_q \cdot m_j = \left( \frac {\alpha_q}{\beta_q}\right)^{j-1} m_j, \quad j = 0, \, 1, \, 2. \]

\noindent As $\alpha_q \neq \beta_q$ we must have $\# H^0 \left( G_q, \, \text{ad}^0 \overline \rho \right) = \#k$.

Globally $H^0 \left( G_{\mathbb Q}, \, \text{ad}^0 \overline \rho \right)$ and $H^0 \left( G_{\mathbb Q}, \, \text{ad}^0 \overline \rho(1) \right)$ are trivial since $\overline \rho$ is absolutely irreducible -- which is not true when $\ell = 3$ and $\overline \rho$ is not absolutely irreducible over $\mathbb Q(\sqrt{-3})$ -- and Proposition \ref{goins_lemma} states $\# H^1_f \left( G_\ell, \, \text{ad}^0 \overline \rho \right) = \#k \cdot \# H^0 \left( G_\ell, \, \text{ad}^0 \overline \rho \right)$, so putting all of this together we have the desired formula $\# H^1_\Sigma \left( \mathbb Q, \, \text{ad}^0 \overline \rho \right) =  \prod_{q \in \Sigma} \#k \cdot \# H^1_\Sigma \left( \mathbb Q, \, \text{ad}^0 \overline \rho(1) \right)$.  \end{proof}


\section{Modular Deformation Ring}

In this section we introduce a different type of deformation ring using modular forms.  Our goal is to show under certain hypotheses that this modular deformation ring is the same as the universal deformation ring introduced in the previous section.  For an odd prime $\ell$, we continue to fix an $\ell$-adic representation $\rho: G_{\mathbb Q} \to GL_2(\mathcal O)$ that is continuous, ordinary, and ramified at finitely many primes.

\subsection{$\Lambda$-adic Modular Forms}

Fix a positive integer $\kappa$, a positive integer $N = N_0 \, \ell$ in terms of an integer $N_0$ prime to $\ell$, and a Dirichlet character $\chi: \left( \mathbb Z / N \, \mathbb Z \right)^\times \to \mathbb C^\times$.  A holomorphic function $f: \mathfrak H \to \mathbb C$ on the upper-half plane is called a classical modular form of weight $\kappa$, level $N$, and nebentype $\chi$  if
\[ f \left( \frac {a \, \tau + b} {c \, \tau + d} \right) = \chi(d) \, \left( c \, \tau + d \right)^\kappa \, f(\tau) \quad \text{for} \quad \begin{pmatrix} a & b \\ c & d \end{pmatrix} \in SL_2 (\mathbb Z), \quad c \equiv 0 \, (N). \]

\noindent Such modular forms have a $q$-expansion $f(\tau) = \sum_n a_n \, q^n$ in terms of $q = e^{2 \pi \sqrt{-1} \, \tau}$ and $a_n \in \overline {\mathbb Q} \hookrightarrow \overline {\mathbb Q}_\ell$ for some embedding.  We define the space of classical cusp forms $S_\kappa(N, \chi)$ as those modular forms with coefficients in $\mathcal O$ such that $a_n = 0$ when $n \leq 0$.  (Equivalently, we assume $\mathcal O$ is always large enough to contain the coefficients of the modular forms of a fixed level $N$.)  For more properties, see \cite[Chap. 1]{MR97k:11080}.

Denote $\Lambda = \mathcal O[[ X ]]$ as the power series ring in the variable $X$.  For each positive rational integer $\kappa$, we have a specialization map $\varphi_\kappa: \Lambda \to \mathcal O$ defined by $1+X \mapsto \left( 1 + \ell \right)^{\kappa}$.  (We can actually define this map for $\ell$-adic integers $\kappa \in \mathbb Z_\ell$, but we will not need this for what follows.)  Denote the Teichmueller character $\omega_\ell: \mathbb Z_\ell^\times \to \mathbb Z_\ell^\times$ defined by $\omega_\ell(d) = \lim_{n \to \infty} d^{\ell^n}$; the image is the $(\ell-1)$-th roots of unity since $\left( \mathbb Z_\ell / \ell^{n+1} \, \mathbb Z_\ell \right)^\times$ has order $\ell^n (\ell - 1)$.  Consider a collection of forms $\sum_n a_n^{(\kappa)} \, q^n \in S_\kappa(N, \chi \, \omega_\ell^{1-\kappa})$ for each $\kappa$.  We call such a collection a Hida family if for each $n$ there exist power series $a_n(X) \in \Lambda$ such that $a_n^{(\kappa)} = a_n \left( (1+\ell)^\kappa - 1 \right)$ for all but finitely many $\kappa$.  The collection $S(N, \chi)$ of formal series $F(X; \, \tau) = \sum_n a_n(X) \, q^n$ are called $\Lambda$-adic cusp forms of level $N$ and nebentype $\chi$ if $F \left( (1+\ell)^\kappa - 1; \, \tau \right) \in S_\kappa \left( N, \chi \, \omega_\ell^{1-\kappa} \right)$ for all but finitely many $\kappa$.  For all $\kappa$, the specializations $\varphi_\kappa \circ F$ are called $\ell$-adic modular forms -- even when they are not classical.  For more properties see \cite[Chap. 7]{MR94j:11044}, although we have altered the notation regarding the nebentype.

\subsection{Hecke Algebras}

Given $f \in S_\kappa(N, \chi \, \omega_\ell^{1-\kappa})$ define the endomorpisms
\[ \begin{aligned} \left( f | T_p \right) (\tau) & = \frac 1p \, \sum_{i=0}^{p-1} f \left( \frac {\tau + i}p \right) + \chi(p) \, \langle p \rangle^{\kappa-1} \, f \left( p \, \tau \right) & \quad & \text{for $p \nmid N$ and} \\ \left( f | U_p \right) (\tau) & = \frac 1p \, \sum_{i=0}^{p-1} f \left( \frac {\tau + i}p \right) & \quad & \text{for $p \mid N$,} \end{aligned} \]

\noindent  in terms of the character $\langle \cdot \rangle: \mathbb Z_\ell^\times \to 1 + \ell \, \mathbb Z_\ell$ mapping $d \mapsto d / \omega_\ell(d)$; and define the $\ell$-adic Hecke algebra $h_\kappa(N, \chi \, \omega_\ell^{1-\kappa})$ as the $\mathcal O$-algebra generated by these operators.  Following \cite[pg. 209]{MR94j:11044}, we extend this action to $\Lambda$-adic modular forms $F(X; \tau) = \sum_n a_n(X) \, q^n$ by considering the $q$-series expansions:
\[ \begin{aligned} \left( F | T_p \right) (X; \, \tau) & = \sum_{p \mid n} a_n(X) \, q^{n/p} + \chi(p) \, \sigma_\ell(p) \,  \sum_n a_n(X) \, q^{pn} & \quad & \text{for $p \nmid N$ and} \\ \left( F | U_p \right) (X; \, \tau) & = \sum_{p \mid n} a_n(X) \, q^{n/p} & \quad & \text{for $p \mid N$,} \end{aligned} \]

\noindent in terms of the character $\sigma_\ell: \mathbb Z_\ell^\times \to \Lambda^\times$ mapping $d \mapsto \left( (1+X)/(1+\ell) \right)^{s(d)}$ where $s(d) = \log \langle d \rangle / \log (1 + \ell) \in \mathbb Z_\ell$; and define the $\Lambda$-adic Hecke algebra $h(N, \chi)$ as the $\Lambda$-algebra generated by these operators.  One checks that the composition $\varphi_\kappa \circ \sigma_\ell$ sends $p \mapsto \langle p \rangle^{\kappa-1}$ for $p \nmid N$, so $\text{ker} \, \varphi_\kappa = \left( \sigma_\ell (1 + \ell) - (1 + \ell)^{\kappa - 1} \right) \Lambda$ is a prime and for any operator $T$ and $\Lambda$-adic cusp form $F$ we have $(\varphi_\kappa \circ F) | T = \varphi_\kappa \circ (F|T)$.

The full $\Lambda$-adic Hecke algebra is too big for our purposes, so following \cite[pg. 202]{MR94j:11044}, we consider a subalgebra.   Set $e = \lim_{n \to \infty} {U_\ell}^{n!} \in h(N, \chi)$ as an idempotent and define the projection $h^0(N, \chi) = e \cdot h(N, \chi)$ as the ordinary part of the Hecke algebra.  Say $F$ is an eigenform with $F | U_\ell = a_\ell(X) \, F$; either $a_\ell(X) \in \Lambda^\times$ is a unit so that $F | e = F$ or $a_\ell(X)$ is in the maximal ideal so $F | e = 0$.  We define ordinary cusp forms as follows: The pairing $S(N, \chi) \times h(N, \chi) \to \Lambda$ defined by mapping a modular form $F$ and an endomorphism $T$ to the first term in the $q$-expansion of $F | T$ induces an $\Lambda$-linear functional $h(N, \chi) \to \Lambda$ defined by $\pi_F = \langle F, \, \cdot \rangle$.  In particular if $F$ is a normalized eigenform for $h(N, \chi)$ i.e. $a_1(X) = 1$ then $\pi_F(T_p) = a_p(X)$ for all $p \nmid N$.  ($\chi$ must be odd i.e. $\chi(-1) = -1$ for normalized eigenforms to exist.)  We now define the ordinary $\Lambda$-adic cusp forms $S^0 \left( N, \chi \right)$ as those cusp forms which make the following diagram commute for all but finitely many $\kappa$:
\[ \begin{CD}  S^0 \left( N, \, \chi \right) @>{\sim}>> \text{Hom}_{\Lambda} \left( h^0(N, \, \chi), \, \Lambda \right) \\ @VVV @VVV \\ S \left( N, \, \chi \right) @>{\sim}>> \text{Hom}_{\Lambda} \left( h(N, \, \chi), \, \Lambda \right) \\ @VV{\varphi_\kappa}V @VV{\varphi_\kappa}V \\ S_\kappa \left( N, \, \chi \, \omega_\ell^{1-\kappa} \right) @>{\sim}>> \text{Hom}_{\mathcal O} \left( h_\kappa(N, \, \chi \, \omega_\ell^{1-\kappa}), \, \mathcal O \right) \end{CD} \]

\subsection{Modular Galois Representations}

We wish to associate Galois representations to modular forms.  To this end, for $\chi$ an odd Dirichlet character modulo $N$ we identify the composition
\[ \begin{CD} G_{\mathbb Q} @>>> \text{Gal} \left( \mathbb Q(\zeta_N) / \mathbb Q \right) @>{\sim}>> \left( \mathbb Z / N \, \mathbb Z \right)^\times @>{\chi}>> \mathcal O^\times \end{CD} \]

\noindent as the Galois representation associated to $\chi$, and identify the compositions $\omega_\ell \circ \varepsilon_\ell$ and $\sigma_\ell \circ \varepsilon_\ell$ in terms of the cyclotomic character as the Galois representations associated to $\omega_\ell$ and $\sigma_\ell$, respectively.  We continue to denote these Galois representations by $\chi$, $\omega_\ell$, and $\sigma_\ell$.  Note that $\omega_\ell$ and $\sigma_\ell$ are unramified at all places $\nu \neq \ell, \, \infty$; and we have the relations $\omega_\ell \equiv \varepsilon_\ell \, (\ell)$ and $\varphi_\kappa \circ \sigma_\ell = (\varepsilon_\ell / \omega_\ell)^{\kappa - 1}$.

Associated to each ordinary normalized classical eigenform $f(\tau) = \sum_n a_n \, q^n$ of weight $\kappa$, level $N$ and nebentype $\chi \, \omega_\ell^{1-\kappa}$ there is a continuous $\ell$-adic Galois representation $\rho_f: G_{\mathbb Q} \to GL_2 \left( \overline {\mathbb Q}_\ell \right)$ with the properties
\begin{enumerate}
\item $\rho_f$ is unramified outside of the primes that divide $N$;
\item $\text{tr} \, \rho_f \left( \text{Frob}_p \right) = a_p$ for $p \nmid N$; and
\item $\det \rho_f = \chi \cdot \left( \varepsilon_\ell / \omega_\ell \right)^{\kappa - 1}$; and
\item  $\rho_f$ is ordinary.
\end{enumerate}

\noindent (Recall that $\ell$ divides $N$ by assumption.) This is due to Eichler and Shimura \cite[Chap. 7]{MR47:3318} when $\kappa = 2$; Deligne and Serre \cite{MR52:284} when $\kappa = 1$; and Deligne \cite{MR42:7460} for $\kappa > 2$.  For proofs, see \cite{MR57:16310}.  Similarly, associated to each ordinary normalized $\Lambda$-adic eigenform $F(X; \tau) = \sum_n a_n(X) \, q^n$ of level $N$ and nebentype $\chi$ there is a continuous $\Lambda$-adic Galois representation $\rho_F: G_{\mathbb Q} \to GL_2 \left( \overline {\mathbb Q}_\ell[[X]] \right)$ with the properties
\begin{enumerate}
\item $\rho_F$ is unramified outside of the primes that divide $N$;
\item $\text{tr} \, \rho_F \left( \text{Frob}_p \right) = a_p(X)$ for $p \nmid N$;
\item $\det \rho_F = \chi \cdot \sigma_\ell$; and
\item $\rho_F$ is ordinary.
\end{enumerate}

\noindent This is due to Hida \cite[Theorem 2.1]{MR87k:11049}; see also \cite[\S 7.5, Theorem 3]{MR94j:11044}.  As $F$ specializes to $\ell$-adic cusp forms $f = \varphi_\kappa \circ F$ the Galois representation $\rho_F$ specializes to $\ell$-adic representations $\rho_f = \varphi_\kappa \circ \rho_F$.

Fix a finite set $\Sigma$ of places that does not contain $\ell$.  Let us consider an $\ell$-adic representation $\rho: G_{\mathbb Q} \to GL_2 \left( \mathcal O \right)$ that is ordinary and ramified at finitely many primes where $\overline \rho$ is odd, absolutely irreducible and wildly ramified at $\ell$.  Fix $\kappa$ and $\chi$ by the relation $\chi \cdot \left( \varepsilon_\ell / \omega_\ell \right)^{\kappa - 1} = \det \rho$.  We say that an ordinary normalized $\Lambda$-adic eigenform $F$ is a cusp form for $\overline \rho$ of type $\Sigma$ if $\varphi_\kappa \circ \rho_F$ is a deformation of $\overline \rho$ of type $\Sigma$.  The level $N_F$ of such a cusp form is divisible by $N_\emptyset = N(\overline \rho) \, \ell$ yet divides $N_\Sigma = N_\emptyset \, \prod_p p^2$ in terms of the product over $p \in \Sigma$ not dividing $N_\emptyset$.   The assumption that $\overline \rho$ is wildly ramified at $\ell$ forces $\ell$ to divide the level $N_\emptyset$; this manifest in the concept of companion forms.  For more properties, see \cite{MR91i:11060} or \cite{MR93i:11063}.

Define a map $\pi: h^0(N_\Sigma, \chi) \to \prod_{\text{$F$ type $\Sigma$}} \Lambda$ by
\[ \pi: \quad T_p \mapsto \left( \dots, \, \pi_F(T_p) , \dots \right) = \left( \dots, \, \text{tr} \, \rho_F \left( \text{Frob}_p \right) , \dots \right) \quad \text{for} \quad p \nmid N_\Sigma; \]

\noindent where each component corresponds to a cusp form for $\overline \rho$ of type $\Sigma$.  We define the modular deformation ring $\mathbb T_\Sigma$ to be the $\Lambda$-algebra generated by the images of $T_p$ for $p \nmid N_\Sigma$.  Note there is a continuous representation
\[ \rho^{\text{mod}}_\Sigma: G_{\mathbb Q} \to GL_2 \left( \mathbb T_\Sigma \right) \]

\noindent where $\rho^{\text{mod}}_\Sigma \simeq \prod_{\text{$F$ type $\Sigma$}} \rho_F$ is a deformation of $\overline \rho$ of type $\Sigma$.

\subsection{Modularity Criteria}

In order to apply the results of the last section, we list a few properties of the modular deformation ring $\mathbb T_\Sigma$.

\begin{localization} \label{localization}  For $\ell$ an odd prime, let $\rho: G_{\mathbb Q} \to GL_2 \left( \mathcal O \right)$ be a continuous $\ell$-adic representation such that $\rho$ is ordinary and ramified at finitely many primes, while $\overline \rho$ is modular, absolutely irreducible, and wildly ramified at $\ell$.  Fix a finite set $\Sigma$ of places that does not contain $\ell$.
\begin{enumerate} 
\item There exist $\Lambda$-adic cusp forms for $\overline \rho$ of type $\Sigma$.
\item $\mathbb T_\Sigma \simeq h^0 \left( N_\Sigma, \, \chi \right)_{\mathfrak m_\Sigma}$ is the localization at some maximal ideal $\mathfrak m_\Sigma$.  
\item $\mathbb T_\Sigma$ is a finitely generated, torsion-free, local $\Lambda$-algebra.  
\end{enumerate} \end{localization}

\begin{proof}  We prove (1).  With $\kappa$, $\chi$, and $N_\emptyset$ as above, fix an integer $\kappa(\overline \rho)$ satisfying $2 \leq \kappa (\overline \rho) \leq \ell$ and $\kappa(\overline \rho) \equiv \kappa \mod (\ell-1)$.  By the ``level lowering'' arguments in \cite{MR91g:11066} and \cite{MR95d:11056} we can find $f \in S_{\kappa(\overline \rho)} \left( N_\emptyset, \, \chi \, \omega_\ell^{1-\kappa(\overline \rho)} \right)$ such that $\overline \rho \simeq \overline \rho_f$.  Choose a $\Lambda$-adic cusp form $F \in S^0(N_\emptyset, \, \chi)$ such that $f = \varphi_{\kappa(\overline \rho)} \circ F$.  Then $F$ is a cusp form of type $\emptyset$, and hence of type $\Sigma$.

We show how (3) follows from (2).  The localization $h^0 \left( N_\Sigma, \, \chi \right)_{\mathfrak m_\Sigma}$ is finitely generated since $h^0 \left( N_\Sigma, \, \chi \right)$ is finite-dimensional.  It is well-known that $h^0 \left( N_\Sigma, \, \chi \right)$ is torsion-free; see \cite[Theorem 3.1]{MR88i:11023} or \cite[\S 7.3, Theorem 1]{MR94j:11044}.  Hence the statements follow for $\mathbb T_\Sigma \simeq h^0 \left( N_\Sigma, \, \chi \right)_{\mathfrak m_\Sigma}$.

We now prove (2).  First we express the $\Lambda$-adic Hecke algebra as a ``complete intersection.''  Denote $\mathcal L$ as the field of fractions of $\Lambda$.  Fix a cusp form $F(X; \tau) = \sum_n a_n(X) \, q^n \in S^0( N_\Sigma, \chi)$, but say that it is a newform of level $N_F$.   It is well-known that there is a decomposition $S^0 \left(N_\Sigma, \, \chi \right) \otimes_{\Lambda} \mathcal L = \bigoplus_F S_F$ over such newforms, where $S_F$ is the $\mathcal L$-linear span of $\left \{ \left. F( X; \, d \, \tau) \, \right | \, \text{$d$ divides $N_\Sigma / N_F$} \right \}$ of dimension
\[ d_F = \# \left \{ d \, \left | \, \text{$d$ divides $\frac {N_\Sigma}{N_F}$} \right. \right \} = \prod_{p | N_\Sigma / N_F} (m_p+1) \quad \text{with} \quad \frac {N_\Sigma}{N_F} = \prod_{p | N_\Sigma / N_F} p^{m_p}. \]

\noindent Let $A_F = \mathcal L \left[ t_p : \text{$p$ divides $N_\Sigma / N_F$} \right]$ be a polynomial ring over $\mathcal L$, and consider the map $A_F \to \text{End}_{\mathcal L} \, S_F$ defined by $t_p \mapsto T_p$.  For each divisor $d_0$ of $N_\Sigma / N_F$ prime to $p$, one checks that the characteristic polynomial $T_p$ on the $(m_p+1)$-dimensional space $\left \{ \left. F( X; \, d_0 \, p^m \, \tau) \, \right | \, m \leq m_p \right \}$ is $P_p(t) = t^{m_p+1} - a_p(X) \, t^{m_p} + (\chi \cdot \sigma_\ell)(p) \, t^{m_p-1}$, setting $\chi(p) = 0$ if $p$ divides $N_F$.  Hence the kernel contains the ideal $I_F = \left( P_p(t_p) : \text{$p$ divides $N_\Sigma / N_F$} \right)$.  But $\dim_{\mathcal L} A_F / I_F = d_F = \dim_{\mathcal L} \text{End}_{\mathcal L} \, S_F$ so
\[ \begin{CD} \pi: \ h^0 \left( N_\Sigma, \, \chi \right) \otimes_{\Lambda} \mathcal L @>{\sim}>> \prod_F \text{End}_{\mathcal L} \, S_F @>{\sim}>> \prod_F A_F / I_F \end{CD} \]

\noindent  where $\pi =  \prod_F \pi_F$ has components $\pi_F$ (corresponding to each newform) mapping $T_p \mapsto t_p \mod I_F$ if $p$ divides $N_\Sigma / N_F$ and $T_p \mapsto a_p(X)$ otherwise.

Next we show $\mathbb T_\Sigma \otimes_{\Lambda} \mathcal L$ is a localization of the $\Lambda$-adic Hecke algebra.  For each maximal ideal $\mathfrak m$ in $h^0(N_\Sigma, \chi)$ there is a natural isomorphism of localizations 
\begin{equation} \label{hecke_localization} h^0 \left( N, \, \chi \right)_{\mathfrak m} \otimes_{\Lambda} \mathcal L \simeq \prod_{\mathfrak p} h^0 \left( N, \, \chi \right)_{\mathfrak p} \simeq \prod_F \prod_{J_F} \left( A_F / I_F \right)_{J_F / I_F} \end{equation}

\noindent where the product is over newforms $F$ and primes $J_F / I_F \subseteq A_F / I_F$ corresponding to primes $\mathfrak p = \pi_F^{-1} (J_F / I_F) \subseteq h^0 (N_\Sigma, \chi ) \otimes_{\Lambda} \mathcal L$ such that $\mathfrak p \cap h^0 ( N_\Sigma, \chi ) \subseteq \mathfrak m$.  Both $\rho$ and $\rho_F$ are unramified at $p$ for $p \nmid N_\Sigma$, so Tchebotarev density implies $F$ is a cusp form for $\overline \rho$ of type $\Sigma$ if and only if the composition
\[ \begin{CD} h^0 \left( N, \, \chi \right) @>{\pi_F}>> \Lambda @>{\varphi_\kappa}>> \mathcal O @>{\mod \lambda}>> k \end{CD} \]

\noindent maps $T_p \mapsto \text{tr} \, \overline \rho \left( \text{Frob}_p \right)$ for all $p \nmid N_\Sigma$.  (Recall that $\kappa$ is the weight associated to $\rho$.)   Choose $\mathfrak m_\Sigma$ as the maximal ideal of $h^0(N_\Sigma, \chi)$ containing $\text{ker} \, \pi_F$ for such cusp forms.  Then the product in \eqref{hecke_localization} is actually over newforms $F$ of type $\Sigma$, where $J_F = \left( t_p : \text{$p$ divides $N_\Sigma / N_F$} \right) \supseteq I_F$.  But then $A_F / J_F \simeq \mathcal L$ implies $\left( A_F / I_F \right)_{J_F / I_F} \simeq \mathcal L$ so we have the isomorphism
\[ h^0 \left( N, \, \chi \right)_{\mathfrak m_\Sigma} \otimes_{\Lambda} \mathcal L \simeq \prod_{\text{$F$ type $\Sigma$}} \left( A_F / I_F \right)_{J_F / I_F} \simeq \prod_{\text{$F$ type $\Sigma$}} \mathcal L \simeq \mathbb T_\Sigma \otimes_{\Lambda} \mathcal L; \]

\noindent this time as the product over cusp forms for $\overline \rho$ of type $\Sigma$.

Finally we show $h^0(N, \chi)_{\mathfrak m_\Sigma} \hookrightarrow \mathbb T_\Sigma$.  For $p \nmid N_\Sigma$, each $T_p \in h^0(N, \chi)_{\mathfrak m_\Sigma}$ maps to $\pi(T_p) \in \mathbb T_\Sigma \cap \prod_{\text{$F$ type $\Sigma$}} (A_F / I_F)_{J_F / I_F}$, so it suffices to show the same for $p$ dividing $N_\Sigma$.  Fix a cusp form $F(X; \tau) = \sum_n a_N(X) \, q^n$ of type $\Sigma$.  If $p$ divides $N_\Sigma / N_F$ then $T_p \mapsto 0$ so consider $p$ that divides $N_F$ but not $N_\Sigma / N_F$.  The projection $\mathbb T_\Sigma \to \Lambda$ gives the composition
\[ \begin{CD} \rho_F : \quad G_{\mathbb Q} @>{\rho^{\text{mod}}_\Sigma}>> GL_2 \left( \mathbb T_\Sigma \right) @>>> GL_2 (\mathcal L) \end{CD} \]

\noindent which shows $\pi_F( T_p) = a_p(X) = \text{tr} \, \left( \rho_F \right)^{I_p} (\text{Frob}_p) =  \text{tr} \, \left( \rho^{\text{mod}}_\Sigma \right)^{I_p} (\text{Frob}_p)$ is an element of $\mathbb T_\Sigma \cap \prod_{\text{$F$ type $\Sigma$}} (A_F / I_F)_{J_F / I_F}$. Hence $h^0(N, \chi)_{\mathfrak m_\Sigma} \simeq \mathbb T_\Sigma$.  \end{proof}

The modular deformation ring $\mathbb T_\Sigma$ is a complete, Noetherian, local $\Lambda$-algebra so there is a unique $\Lambda$-algebra surjection $\phi_\Sigma: R_\Sigma \to \mathbb T_\Sigma$ of the universal deformation ring such that $\rho^{\text{mod}}_\Sigma \simeq \phi_\Sigma \circ \rho^{\text{univ}}_\Sigma$.  The following result states that this map is an isomorphism if and only if it is an isomorphism upon specializing the weight.

\begin{isomorphism}  \label{isomorphism} For $\ell$ an odd prime, let $\rho: G_{\mathbb Q} \to GL_2(\mathcal O)$ be a continuous $\ell$-adic representation such that $\rho$ is ordinary and ramified at finitely many primes, while $\overline \rho$ is absolutely irreducible, modular, and wildly ramified at $\ell$.  For each positive integer $\kappa$, denote $\mathbb T^{(\kappa)}_\Sigma$ and $R^{(\kappa)}_\Sigma$ as the images of $\mathbb T_\Sigma$ and $R_\Sigma$ under the specialization map $\varphi_\kappa$, respectively.  Then the following are equivalent:
\begin{enumerate}
\item $\phi_\Sigma: R_\Sigma \to \mathbb T_\Sigma$ is an isomorphism.
\item $\phi^{(\kappa)}_\Sigma: R^{(\kappa)}_\Sigma \to \mathbb T^{(\kappa)}_\Sigma$ is an isomorphism for all positive integers $\kappa$.
\item $\phi^{(\kappa)}_\Sigma: R^{(\kappa)}_\Sigma \to \mathbb T^{(\kappa)}_\Sigma$ is an isomorphism for some positive integer $\kappa$.
\end{enumerate}

\noindent If any of the above hold for every finite set $\Sigma$ not containing $\ell$, then $\rho$ is $\ell$-adically modular i.e. $\rho \simeq \rho_f$ for some $\ell$-adic cusp form $f$.  \end{isomorphism}

\begin{proof} We prove the last statement of the proposition assuming $\phi_\Sigma: R_\Sigma \to \mathbb T_\Sigma$ is an isomorphism for every finite set $\Sigma$.  Choose $\Sigma$ large enough so that $\rho$ is a deformation of $\overline \rho$ of type $\Sigma$, define $\kappa$ and $\chi$ by the relation $\det \rho = \chi \cdot \left( \varepsilon_\ell / \omega_\ell \right)^{\kappa - 1}$.  There is a surjection $\phi: R_\Sigma \to \Lambda$, so there is a map $\pi_\Sigma : \mathbb T_\Sigma \to \Lambda$ such that $\phi = \pi_\Sigma \circ \phi_\Sigma$.  However, the maps $\mathbb T_\Sigma \to \Lambda$ all correspond to the canonical pairing $S(N_\emptyset, \chi) \times h(N_\emptyset, \chi) \to \Lambda$ so $\pi_\Sigma \simeq \pi_F$ for some cusp form $F$ of type $\Sigma$.  Then $f = \varphi_\kappa \circ F$ is the desired $\ell$-adic cusp form because
\[ \rho \simeq \left( \varphi_\kappa \circ \phi \right) \circ \rho^{\text{univ}}_\Sigma \simeq  \varphi_\kappa \circ \pi_F \circ  \rho^{\text{mod}}_\Sigma \simeq \varphi_\kappa \circ \rho_F \simeq \rho_f. \]

Now we prove the equivalence of the statements; we are motivated by the proof of \cite[Proposition 3.4]{MR1845181}.  For each positive integer $\kappa$ we have the commutative diagram
\[ \begin{CD} 0 @>>> \text{ker} \, \varphi_\kappa \cdot R_\Sigma @>>> R_\Sigma @>{\varphi_\kappa}>> R^{(\kappa)}_\Sigma @>>> 0 \\ @. @VV{\phi_\Sigma}V @VV{\phi_\Sigma}V @VV{\phi^{(\kappa)}_\Sigma}V  \\ 0 @>>> \text{ker} \, \varphi_\kappa \cdot \mathbb T_\Sigma @>>> \mathbb T_\Sigma @>{\varphi_\kappa}>> \mathbb T^{(\kappa)}_\Sigma @>>> 0 \end{CD} \]

\noindent where $\text{ker} \, \varphi_\kappa = \left( \sigma_\ell(1+\ell) - (1+\ell)^{\kappa - 1} \right) \Lambda$ is a prime ideal.  Clearly from the diagram (1) implies (2), and (2) obviously implies (3).  Assume (3) holds i.e. $\phi^{(\kappa)}_\Sigma$ is an isomorphism for some positive integer $\kappa$.  Then $R_\Sigma / \text{ker} \, \varphi_\kappa \cdot R_\Sigma \simeq \mathbb T_\Sigma / \text{ker} \, \varphi_\kappa \cdot T_\Sigma$ as integral domains, but by Proposition \ref{localization} both $R_\Sigma$ and $\mathbb T_\Sigma$ are torsion-free so $R_\Sigma \simeq \mathbb T_\Sigma$.  \end{proof}

\subsection{Specialization to Weight 2}

We now show that under suitable hypotheses the modular deformation ring is the same as the universal deformation ring.  This is the main result of the paper.

\begin{wiles} \label{wiles} For $\ell$ an odd prime, let $\rho: G_{\mathbb Q} \to GL_2(\mathcal O)$ be a continuous $\ell$-adic representation such that
\begin{enumerate}
\item $\rho$ is ordinary and ramified at finitely many primes;
\item $\overline \rho$ is absolutely irreducible when restricted to $\text{Gal} \left( \overline {\mathbb Q} / \mathbb Q(\sqrt{(-1)^{(\ell-1)/2} \, \ell}) \right)$, modular, and wildly ramified at $\ell$.
\end{enumerate}

\noindent Then $\rho$ is $\ell$-adically modular i.e. $\rho \simeq \rho_f$ for an $\ell$-adic cusp form $f$.
\end{wiles}

\begin{proof} By Proposition \ref{isomorphism} it suffices to show $\phi^{(2)}_\Sigma: R^{(2)}_\Sigma \to \mathbb T^{(2)}_\Sigma$ (i.e. the weight $\kappa = 2$ case) is an isomorphism for all finite sets $\Sigma$ not containing $\ell$.  We identify $\mathbb T^{(2)}_\Sigma$ as the modular deformation ring associated to ordinary cusp forms $f(\tau) = \sum_n a_n \, q^n$ of weight 2, level $N_\Sigma$, and nebentype $\chi \cdot \omega_\ell^{-1}$ such that $\overline \rho \simeq \overline \rho_f$; and identify $R^{(2)}_\Sigma$ is its universal deformation ring.  Recall that $N_\Sigma = N_\emptyset \, \prod_p p^2$ in terms of the product over $p \in \Sigma$ not dividing $N_\emptyset = N(\overline \rho) \, \ell$.  The following commutative diagram is exact for the unique surjection $\phi$:
\[ \begin{CD} 0 @>>> \text{ker} \, \phi @>>> R^{(2)}_\Sigma @>{\phi}>> \mathcal O @>>> 0 \\ @. @VVV @VV{\phi^{(2)}_\Sigma}V @VVV \\ 0 @>>> \text{ker} \, \pi_f @>>> \mathbb T^{(2)}_\Sigma @>{\pi_f}>> \mathcal O @>>> 0 \end{CD} \]

We fix some notation.  Choose a cusp form $f(\tau)$ of level $N_\emptyset$.  We will denote the ``tangent space'' of $R^{(2)}_\Sigma$ as $\Phi_\Sigma = \left( \text{ker} \, \phi \right) / \left( \text{ker} \, \phi \right)^2$, as well as the ideals $\mathfrak p_\Sigma = \text{ker} \, \pi_f$ and $I_\Sigma = \text{Ann}_{\mathbb T^{(2)}_\Sigma} \text{ker} \, \pi_f$ in $\mathbb T^{(2)}_\Sigma$.  We assume for the moment that there exists a family of $\mathbb T^{(2)}_\Sigma$-modules $L_\Sigma$ satisfying the following properties:
\begin{enumerate}
\item[HM1:] $L_\Sigma$ is free over $\mathcal O$ where $\text{rank}_{\mathcal O} L_\Sigma = 2 \cdot \text{rank}_{\mathcal O} \mathbb T^{(2)}_\Sigma$ and $\text{rank}_{\mathcal O} L_\Sigma[\mathfrak p_\Sigma] = 2$.
\item[HM2:] There is a pairing $\langle \cdot, \cdot \rangle_\Sigma: L_\Sigma \times L_\Sigma \to \mathcal O$ such that $L_\Sigma \simeq \text{Hom}_{\mathcal O} \left( L_\Sigma, \mathcal O \right)$.
\item[HM3:] For $\Sigma' \subseteq \Sigma$, there exist surjective $\mathbb T^{(2)}_\Sigma$-homomorphisms $\beta_{\Sigma'} : L_\Sigma \to L_{\Sigma'}$.
\end{enumerate}

\noindent The proof of this theorem will be divided into three steps: first we show how the minimal case follows from the existence of these Hecke modules, second we show how the general case follows from the minimal case, and third we construct these Hecke modules using modular curves.  The arguments below follow the discussion in \cite{MR98c:11047} and \cite{MR1639612}, where we use an axiomatic approach following \cite{MR2002m:11042}.

\emph{The Minimal Case.}  We show $\phi^{(2)}_\emptyset: R^{(2)}_\emptyset \to \mathbb T^{(2)}_\emptyset$ is an isomorphism.  Denote $\overline {R}^{(2)}_\emptyset = R^{(2)}_\emptyset / \lambda \, R^{(2)}_\emptyset$ and $\overline {\mathbb T}^{(2)}_\emptyset = \mathbb T^{(2)}_\emptyset / \lambda \, \mathbb T^{(2)}_\emptyset$ as $k$-vector spaces.  It is well-known that $R^{(2)}_\emptyset \to \mathbb T^{(2)}_\emptyset$ is an isomorphism of complete intersections if and only if $\overline {R}^{(2)}_\emptyset \to \overline {\mathbb T}^{(2)}_\emptyset $ is an isomorphism of complete intersections; see \cite[Lemma 5.29]{MR99d:11067b}.  We show the criteria of  \cite[Theorem 2.1]{MR98c:11047} are satisfied.

We construct a sequence of surjective maps.  Denote $r = \dim_k H^1_\emptyset \left( \mathbb Q, \, \text{ad}^0 \overline \rho(1) \right)$.  For each positive integer $n$, let $Q$ be a set of $r$ primes $q$ satisfying
\begin{enumerate}
\item $q \equiv 1 \mod{\ell^n}$;
\item $\overline \rho$ is unramified at $q$ and $\overline \rho \left( \text{Frob}_q \right)$ has distinct $k$-rational eigenvalues; and
\item $H^1_\emptyset \left( \mathbb Q, \, \text{ad}^0 \overline \rho \right)$ may be embedded in $\bigoplus_{q \in Q} H^1 \left( G_q, \, \text{ad}^0 \overline \rho \right)$.
\end{enumerate}

\noindent Such a collection $Q$ exists since we assume $\overline \rho$ remains absolutely irreducible over $\mathbb Q ( \sqrt{(-1)^{(\ell-1)/2} \, \ell} )$; see \cite[Theorem 2.49]{MR99d:11067b} and its proof.  By Proposition \ref{dimension} -- note $\rho$ is odd because it is residually modular -- we have $\dim_k H^1_Q \left( \mathbb Q, \, \text{ad}^0 \overline \rho \right) = r + \dim_k H^1_Q \left( \mathbb Q, \, \text{ad}^0 \overline \rho(1) \right)$.  But 
\[ H^1_Q \left( \mathbb Q, \, \text{ad}^0 \overline \rho(1) \right) = \text{ker} \left[ H^1_\emptyset \left( \mathbb Q, \, \text{ad}^0 \overline \rho \right) \to \bigoplus_{q \in Q} H^1 \left( G_q, \, \text{ad}^0 \overline \rho \right) \right] \]

\noindent is trivial so $\dim_k H^1_Q \left( \mathbb Q, \, \text{ad}^0 \overline \rho \right) = r$.  Using Proposition \ref{dimension} again we see $R^{(2)}_Q$ may be generated by $r$ elements as an $\mathcal O$-algebra.  The universal property implies that there is a surjection $R^{(2)}_Q \to R^{(2)}_\emptyset$, so $\overline {R}^{(2)}_\emptyset$ has dimension at most $r$ as a $k$-vector space.  Denote $\Delta_Q$ as the maximal quotient of $\text{Gal} \left( \prod_{q \in Q} \mathbb Q(\zeta_q)/ \mathbb Q \right) \simeq \prod_{q \in Q} \left( \mathbb Z / q \, \mathbb Z \right)^\times$ of $\ell$-power order.  Given $r$ variables $X_i$ we have a surjection of the power series ring $\mathcal O[[X_1, \dots, \, X_r]] \to \overline {R}^{(2)}_\emptyset$.  Introduce a second power series ring $\mathcal O[[S_1, \dots, \, S_r]]$ with prime ideal $\mathfrak a = \left( S_1, \dots, \, S_r \right)$, and a map onto $\mathcal O[\Delta_Q]$ that sends $1+S_i$ to one of its $r$ generators; the kernel of this map is contained in $\mathfrak a^n$.  The following diagram commutes:
\[ \begin{CD} \mathcal O[[S_1, \dots, \, S_r]] @>>> \mathcal O[\Delta_Q] \\ @VVV @VVV \\ \mathcal O[[X_1, \dots, \, X_r]] @>>> \overline {R}^{(2)}_\emptyset @>{\phi^{(2)}_\emptyset}>> \overline {\mathbb T}^{(2)}_\emptyset \end{CD} \]

\noindent where the horizontal maps are surjections, and the vertical maps are chosen so that the image of $\mathfrak a$ in $R^{(2)}_\emptyset$ is trivial.

Now impose the following additional condition on the Hecke modules:
\begin{enumerate}
\item[HM5:] $L_Q$ is free over $\mathcal O[\Delta_Q]$ and $L_Q / \mathfrak a \, L_Q \simeq L_\emptyset$.
\end{enumerate}

\noindent (The numbering of this property will become clear later.)  Denote $L_Q' = L_Q / \mathfrak a^n \, L_Q$ so that $L_Q' / \mathfrak a \, L_Q' \simeq L_\emptyset$.  Property (HM5) implies $\text{Ann}_{\mathcal O[[S_1, \dots, \, S_r]]} L_Q \subseteq \mathfrak a^n$, so we actually have $\text{Ann}_{\mathcal O[[S_1, \dots, \, S_r]]} L_Q' = \mathfrak a^n$.  We have the surjections $L_Q \to L_Q'$ and $\mathcal O[\Delta_Q] \to \mathcal O[[S_1, \dots, \, S_r]] / \mathfrak a^n$, so $L_Q'$ is free module over $\mathcal O[[S_1, \dots, \, S_r]] / \mathfrak a^n$.  Then by \cite[Theorem 2.1]{MR98c:11047} where we use property (HM1) and tensor with $k$, both $\overline {R}^{(2)}_\emptyset$ and $\overline {\mathbb T}^{(2)}_\emptyset$ are complete intersections while $L_\emptyset \otimes_{\mathcal O} k$ is free over $\overline {R}^{(2)}_\emptyset$.  However $\text{ker} \left[ \overline {R}^{(2)}_\emptyset \to \overline {\mathbb T}^{(2)}_\emptyset \right] \subseteq \text{Ann}_{\overline {R}^{(2)}_\emptyset} L_\emptyset \otimes_{\mathcal O} k = \{ 0 \}$ so $\overline {R}^{(2)}_\emptyset \simeq \overline {\mathbb T}^{(2)}_\emptyset$.

\emph{Reduction to the Minimal Case.}  We now show $\phi^{(2)}_\Sigma: R^{(2)}_\Sigma \to \mathbb T^{(2)}_\Sigma$ is an isomorphism by verifying the criteria of \cite[Theorem 2.4]{MR98c:11047}.  Denote $\text{ad}^0 \rho_{f,n} = \text{ad}^0 \rho_f \otimes_{\mathcal O} \lambda^{-n} \, \mathcal O / \mathcal O$ so that as a generalization to \eqref{mazur},
\[ \# \, \Phi_\Sigma = \lim_{n \to \infty} \# \, \text{Hom}_{\mathcal O} \left( \Phi_\Sigma, \, \lambda^{-n} \, \mathcal O / \mathcal O \right) = \lim_{n \to \infty} \# \, H_\Sigma^1 \left( \mathbb Q, \, \text{ad}^0 \rho_{f,n} \right). \]

\noindent Hence the change in size from $\Phi_\emptyset$ to $\Phi_\Sigma$ satisfies
\[ \frac {\# \Phi_\Sigma}{\# \, \Phi_\emptyset} \leq \lim_{n \to \infty} \prod_{p \in \Sigma} \frac {\# \, H^1 \left( G_p, \, \text{ad}^0 \rho_{f,n} \right)}{\# \, H_f^1 \left( G_p, \, \text{ad}^0 \rho_{f,n} \right)} =  \lim_{n \to \infty} \prod_{p \in \Sigma} \# \, H^0 \left( G_p, \, \text{ad}^0 \rho_{f,n}(1) \right) \]

\noindent because $\# H_f^1 \left( G_p, \, \text{ad}^0 \rho_{f,n} \right) = \# H^0 \left( G_p, \, \text{ad}^0 \rho_{f,n} \right)$ when $p \neq \ell$, and we only consider finite sets $\Sigma$ not containing $\ell$ since $\overline \rho$ is ramified at $\ell$.

We recall the value of $c_p = \lim_{n \to \infty} \# H^0 \left( G_p, \, \text{ad}^0 \rho_{f,n}(1) \right)$.  Denote $K \subseteq \overline {\mathbb Q}_\ell$ as the field of fractions of $\mathcal O$, and fix a normalized valuation $v_\lambda: K^\times \to \mathbb Z$.  Note that for $c_p \neq 0$ we have $v_\lambda(c_p) = v_\lambda(c_p')$ if and only if $c_p \, \mathcal O = c_p' \, \mathcal O$, so it suffices to compute $v_\lambda(c_p)$; we do so via determinants.  If $\overline \rho$ is unramified at $p$ we have
\[ \det \left[ 1_3 - \text{ad}^0 \rho_f(\text{Frob}_p) \cdot X \right] = (1 - X) \left( (1+X)^2 - \frac {a_p^2}{(\chi \cdot \omega_\ell^{-1})(\text{Frob}_p)} \, \frac Xp \right), \]

\noindent so for any prime $p \not \in \Sigma$ such that $p \neq \ell$ define $t_p \in \mathbb T^{(2)}_\Sigma$ as the element 
\[ t_p = \begin{cases} (1-p) \left( (1+p)^2 - T_p^2 / (\chi \cdot \omega_\ell^{-1})(\text{Frob}_p) \right) & \text{if $\dim_k \overline \rho^{I_p} = 2$;} \\ 1-p^2 & \text{if $\dim_k \overline \rho^{I_p} = 1$, $\det \overline \rho \vert_{I_p} = 1$;} \\ 1-p & \text{if $\dim_k \overline \rho^{I_p} = 1$, $\det \overline \rho \vert_{I_p} \neq 1$;} \\ 1+p & \text{if $\dim_k \overline \rho^{I_p} = 0$, $p  \in P$;} \\ 1 & \text{if $\dim_k \overline \rho^{I_p} = 0$, $p \not \in P$;} \end{cases} \]

\noindent where we denote $P$ as the collection of primes $p \equiv -1 \pod{\ell}$ such that $\overline \rho \vert_{G_p}$ is irreducible yet $\overline \rho \vert_{I_p}$ is reducible.  It is straight-forward to verify
\[ v_\lambda(c_p) = v_\lambda \left( \det \left[ 1_3 - \text{ad}^0 \rho_f(1)^{I_p} (\text{Frob}_p) \right] \right) = v_\lambda \left( t_p \, \text{mod} \, \mathfrak p_\Sigma \right). \]

Now impose the following additional condition on the Hecke modules:
\begin{enumerate}
\item[HM4:] $\left( \beta_{\Sigma'} \circ \widehat \beta_{\Sigma'} \right) L_{\Sigma'} = \prod_{p \in \Sigma - \Sigma'} t_p \cdot L_{\Sigma'}$ where $\widehat \beta_{\Sigma'}: L_{\Sigma'} \to L_\Sigma$ is the adjoint.
\end{enumerate}

\noindent To be precise, the adjoint is with respect to the pairings: $\langle x, \widehat \beta_{\Sigma'} x' \rangle_\Sigma = \langle \beta_{\Sigma'} x, x' \rangle_{\Sigma'}$ for all $x \in L_\Sigma$ and $x' \in L_{\Sigma'}$.  We've seen that property (HM5) implies $R^{(2)}_\emptyset$ is a complete intersection and $L_\emptyset$ is a free $R^{(2)}_\emptyset$-module, while property (HM1) implies $\mathfrak p_\Sigma$ and $\mathfrak p_\emptyset$ are in the support of $L_\Sigma$ and $L_\emptyset$, respectively.  By property (HM2) we have the isomorphism $L_\Sigma / L_\Sigma[I_\Sigma] \simeq \text{Hom}_{\mathcal O} \left( L_\Sigma[\mathfrak p_\Sigma], \mathcal O \right)$, so for any basis $\{ x_1, x_2 \}$ of $L_\Sigma[\mathfrak p_\Sigma]$ we have
\[ \Omega_\Sigma = \frac {L_\Sigma}{L_\Sigma[\mathfrak p_\Sigma] + L_\Sigma[I_\Sigma]} \simeq \frac {\text{Hom}_{\mathcal O} \left( L_\Sigma[\mathfrak p_\Sigma], \mathcal O \right)}{L_\Sigma[\mathfrak p_\Sigma]} \simeq \frac {\mathcal O}{\det \left[ \langle x_i, x_j \rangle_\Sigma \right]_{i,j} \, \mathcal O}. \]

\noindent (Compare with the discussion in \cite[\S 4.4]{MR99d:11067b}.)  By properties (HM3) and (HM4), the adjoint $\widehat \beta_{\Sigma'}$ has torsion-free cokernel, so we can choose the basis $\{ x_1', x_2' \}$ of $L_{\Sigma'}[\mathfrak p_{\Sigma'}]$ such that $\{ x_1, x_2 \} = \widehat \beta_{\Sigma'} \{ x_1', x_2' \}$ is a basis for $L_\Sigma[\mathfrak p_\Sigma]$.  This gives
\[ \left \langle x_i, x_j \right \rangle_\Sigma \mathcal O = \left \langle \left( \beta_{\Sigma'} \circ \widehat \beta_{\Sigma'} \right) x_i', x_j' \right \rangle_{\Sigma'} \mathcal O = \prod_{p \in \Sigma - \Sigma'} c_p \cdot \left \langle x_i', x_j' \right \rangle_{\Sigma'} \mathcal O \]

\noindent which implies, for $\Sigma' = \emptyset$, the inequality
\[ \frac {\# \Omega_\Sigma }{\# \, \Omega_\emptyset} = \prod_{p \in \Sigma} \left[ \det \left[ \langle x_i', x_j' \rangle_{\Sigma'} \right]_{i,j} \, \mathcal O : \det \left[ \langle x_i, x_j \rangle_\Sigma \right]_{i,j} \, \mathcal O \right] = \prod_{p \in \Sigma} c_p^2 \geq \frac {\left( \# \Phi_\Sigma \right)^2}{\left( \# \Phi_\emptyset \right)^2}. \]

\noindent By \cite[Theorem 2.4]{MR98c:11047} we find $\# \Omega_\emptyset = \left( \# \Phi_\emptyset \right)^2$ so that we have the inequality $\# \Omega_\Sigma \geq \left( \# \Phi_\Sigma \right)^2$.  Then the same theorem implies $R^{(2)}_\Sigma$ is a complete intersection and $L_\Sigma$ is a free $R^{(2)}_\Sigma$-module.  Hence $\text{ker} \left[ R^{(2)}_\Sigma \to \mathbb T^{(2)}_\Sigma \right] \subseteq \text{Ann}_{R^{(2)}_\Sigma} L_\Sigma = \{ 0 \}$ so $R^{(2)}_\Sigma \simeq \mathbb T^{(2)}_\Sigma$.

\emph{Construction of Hecke Modules.}  To complete the proof, we construct a family of $\mathbb T^{(2)}_\Sigma$-modules $L_\Sigma$, for finite sets $\Sigma$ not including $\ell$,  satisfying the following properties:
\begin{enumerate}
\item[HM1:] $L_\Sigma$ is free over $\mathcal O$ where $\text{rank}_{\mathcal O} L_\Sigma = 2 \cdot \text{rank}_{\mathcal O} \mathbb T^{(2)}_\Sigma$ and $\text{rank}_{\mathcal O} L_\Sigma[\mathfrak p_\Sigma] = 2$.
\item[HM2:] There is a pairing $\langle \cdot, \cdot \rangle_\Sigma: L_\Sigma \times L_\Sigma \to \mathcal O$ such that $L_\Sigma \simeq \text{Hom}_{\mathcal O} \left( L_\Sigma, \mathcal O \right)$.
\item[HM3:] For $\Sigma' \subseteq \Sigma$, there exist surjective $\mathbb T^{(2)}_\Sigma$-homomorphisms $\beta_{\Sigma'} : L_\Sigma \to L_{\Sigma'}$.
\item[HM4:] $\left( \beta_{\Sigma'} \circ \widehat \beta_{\Sigma'} \right) L_{\Sigma'} = \prod_{p \in \Sigma - \Sigma'} t_p \cdot L_{\Sigma'}$ where $\widehat \beta_{\Sigma'}: L_{\Sigma'} \to L_\Sigma$ is the adjoint.
\item[HM5:] $L_Q$ is free over $\mathcal O[\Delta_Q]$ and $L_Q / \mathfrak a \, L_Q \simeq L_\emptyset$.
\end{enumerate}

\noindent We will construct the modules using modular curves.

First we recall a few facts about the modular curve $X(N)$ of level $N = \prod_p p^{e_p}$ in order to motivate a more general definition.  A function $f: \mathfrak H \to \mathbb C$ is a cusp form of weight 2 and level $N$ if and only if $f(\tau) \, d \tau$ is a holomorphic differential on the compact Riemann surface $X(N)$.  The upper-half plane $\mathfrak H \simeq SL_2(\mathbb R) / SO_2(\mathbb R)$, so using the Strong Approximation Theorem \cite[Proposition 3.3.1]{MR97k:11080}, we have the double-coset decomposition
\[ \Gamma(N) \backslash \mathfrak H \simeq \Gamma(N) \backslash SL_2(\mathbb R)  / SO_2(\mathbb R) \simeq GL_2(\mathbb Q) \backslash GL_2(\mathbb A_{\mathbb Q}) / U_N \cdot SO_2(\mathbb R) \, \mathbb R^\times \]

\noindent where $U_N = \prod_p \text{ker} \left[ GL_2(\mathbb Z_p)\to GL_2(\mathbb Z / p^{e_p} \mathbb Z) \right]$.  We may identify $X(N)$ as the compactification of this decomposition by adjoining the cusps.  We may construct a cuspidal automorphic form $\pi: GL_2(\mathbb A_{\mathbb Q}) \to \mathbb C$ of level $U_N$ from a cusp form $f$ of weight 2 and level $N$ by defining
\[ \pi(g) = \frac {a \, d - b \, c}{\left( c \, \tau + d \right)^2} \ f \left( \frac {a \, \tau + b}{c \, \tau + d} \right) \quad \text{where} \quad \tau = \sqrt{-1}, \quad g_{\infty} = \begin{pmatrix} a & b \\ c & d \end{pmatrix} \in GL_2(\mathbb R)^+, \]

\noindent and $g = \gamma \, g_{\infty} \, u$ for $\gamma \in GL_2(\mathbb Q)$ and $u \in U_N \cdot SO_2(\mathbb R) \, \mathbb R^\times$.  In general, given an open, compact subgroup $U \subseteq GL_2(\mathbb A^\infty_{\mathbb Q})$, we define $X_U$ as the compactification of $GL_2(\mathbb Q) \backslash GL_2(\mathbb A_{\mathbb Q}) / U \cdot SO_2(\mathbb R) \, \mathbb R^\times$.  Write $N_\Sigma = \ell \prod_{p \neq \ell} p^{e_p}$, and consider in particular the subgroup
\[ \begin{aligned} U_\Sigma & = \text{ker} \left[ GL_2(\mathbb Z_\ell)\to GL_2(\mathbb F_\ell) \right] \times \prod_{p \in P} \text{ker} \left[ GL_2(\mathbb Z_p) \to GL_2(\mathbb Z / p^{e_p/2} \mathbb Z) \right] \\ & \times \prod_{p \in \Sigma - P} U_p(e_p+ \text{dim}_k \overline \rho^{I_p}) \times \prod_{p \not \in \Sigma \cup P \cup \{ \ell \} } U_p(e_p) \end{aligned} \]

\noindent where we denote
\[ U_p(e) = \left \{ \left. \begin{pmatrix} a & b \\ c & d \end{pmatrix} \in GL_2(\mathbb Z_p) \, \right \vert \, \text{$c \equiv 0, \quad d^{\ell^n} \equiv 1 \pod{p^e}$ for some $n$} \right \}. \]

\noindent We identify cusp forms of weight 2 and level $N_\Sigma$ having coefficients in $\mathcal O$ with those holomorphic differentials in $H^1\left( X_{U_\Sigma}, \mathcal O \right)$.

Second we construct an $\mathcal O$-lattice $V_\Sigma$ with a pairing $( \cdot, \cdot )_\Sigma$ following \cite{MR1639612}.  Define the following subgroup of the finite quotient $GL_2(\mathbb A_{\mathbb Q}^{\infty}) / U_\Sigma$:
\[ G_\Sigma = GL_2(\mathbb F_\ell) \times \prod_{p \in P} \left \{ \left. \begin{pmatrix} a & b \\ c & d \end{pmatrix} \in GL_2(\mathbb Z/ p^{e_p/2} \mathbb Z) \, \right \vert \, \text{$c \equiv 0 \, (p)$ if $p \in \Sigma$} \right \} \]

\noindent and consider absolutely irreducible representations $G_\Sigma \to GL(V_\Sigma)$ for some $\mathcal O$-module $V_\Sigma$; we will choose a specific representation associated to certain multiplicative characters $\chi_p$ of finite order prime to $\ell$ for either $p \in P$ or $p = \ell$.  When $p \in P$ the restriction $\overline \rho \vert_{G_p}$ is induced from a character $\text{Gal} \left( \overline {\mathbb Q}_p / \mathbb Q_{p^2} \right) \to \overline {\mathbb F}_\ell^\times$ defined over the unramified quadratic extension $\mathbb Q_{p^2}$ of $\mathbb Q_p$; see \cite{MR1638490}.  Choose a homomorphism $\chi_p$ such that the restriction of this character to the inertia group at $p$ is the composition 
\[ \begin{CD} I_p @>>> W(\mathbb F_{p^2})^\times @>{\chi_p}>> \mathcal O^\times @>{\text{mod} \, \lambda}>> \overline {\mathbb F}_\ell^\times. \end{CD} \]

\noindent This character induces an irreducible $\phi(p^{e_p/2})$-dimensional representation $\Theta(\chi_p)$ of $GL_2(\mathbb Z / p^{e_p/2} \mathbb Z)$; see \cite[\S 3.2]{MR1639612} for details on the construction.  When $p = \ell$, the representation $\overline \rho$ is ordinary so choose a homomorphism $\chi_\ell$ such that we have the composition
\[ \begin{CD} \det \overline \rho: \quad I_\ell @>{\omega_\ell}>> \mathbb F_\ell^\times @>{\chi_\ell}>> \mathcal O^\times @>{\text{mod} \, \lambda}>> \overline {\mathbb F}_\ell^\times. \end{CD} \]

\noindent If $\chi_\ell = \chi_0$ is the map such that $\chi_0 \mod{\lambda}$ is trivial we have a 1-dimensional representation of $GL_2(\mathbb F_\ell)$ given by the composition $\chi_\ell \circ \det$; otherwise we have an irreducible $(\ell+1)$-dimensional representation $I(\chi_\ell, \chi_0)$.  We choose a $G_\Sigma$-invariant $\mathcal O$-lattice $V_\Sigma$ such that the irreducible representation $G_\Sigma \to GL(V_\Sigma \otimes_{\mathcal O} K)$ is the tensor product over $p \in P$ and $p = \ell$ of the irreducible representations described above.  Enlarging $\mathcal O$ if necessary, we can choose $V_\Sigma$ such that there exists a non-degenerate pairing $( \cdot, \cdot)_\Sigma$ for which $V_\Sigma$ self-dual; see \cite[Lemma 3.3.1]{MR1639612}.

We now define $L_\Sigma$ following the discussion in \cite{MR1639612}.  We recall that by Proposition \ref{localization}, where we specialize to weight $\kappa = 2$, we have $\mathbb T^{(2)}_\Sigma \simeq \mathbb T^{(2)}_{\mathfrak m}$ as the localization of $\mathbb T^{(\kappa)} = \varphi_\kappa \left( h^0 \left( N_\Sigma, \, \chi \right) \right)$ at the maximal ideal $\mathfrak m = \varphi_\kappa \left( \mathfrak m_\Sigma \right)$ containing $\mathfrak p_\Sigma = \text{ker} \, \pi_f$.  We define $L_\Sigma$ as a localization at this maximal ideal:
\[ L_\Sigma = \text{Hom}_{\mathcal O[G_\Sigma]} \left( V_\Sigma , \ H^1 \left( X_{U_\Sigma}, \mathcal O \right) \right)_{\mathfrak m}. \]

\noindent Property (HM1) is essentially \cite[Lemma 5.3.1]{MR1639612}.  Properties (HM2), (HM3), and (HM4) are verified in \cite[Proposition 5.5.1]{MR1639612}. Property (HM5) is verified in \cite[Proposition 5.6.1]{MR1639612}.  This completes the proof of Theorem \ref{wiles}. \end{proof}

\subsection{Specialization to Weight 1}

We have shown that under suitable hypotheses $\rho$ is $\ell$-adically modular, but this cusp form may not be classical.  Our goal in this section is to show when such a form is indeed classical.  We will work with projective $\ell$-adic representations
\[ \begin{CD} \widetilde \rho: \quad G_{\mathbb Q} @>{\rho}>> GL_2 \left( \mathcal O \right) @>>> PGL_2 \left( \mathcal O \right). \end{CD} \]

\begin{buzzard} \label{buzzard} For $\ell$ an odd prime, let $\rho: G_{\mathbb Q} \to GL_2 \left( \mathcal O \right)$ be a continuous Galois representation such that
\begin{enumerate}
\item $\rho$ is ramified at finitely many primes;
\item $\overline \rho$ is absolutely irreducible when restricted to $\text{Gal} \left( \overline {\mathbb Q} / \mathbb Q(\sqrt{(-1)^{(\ell-1)/2} \, \ell}) \right)$, modular, and wildly ramified at $\ell$;
\item $\rho(G_\ell)$ is finite and $\widetilde \rho(G_\ell)$ is cyclic of $\ell$-power order.
\end{enumerate}

\noindent Then $\imath \circ  \rho: G_{\mathbb Q} \to GL_2 \left( \mathbb C \right)$ is modular for each embedding $\imath: K \hookrightarrow \mathbb C$.  \end{buzzard}

The conditions above force $\ell  = 3, \, 5$:  Since $\overline \rho$ is wildly ramified its image is divisible by $\ell$, and since it is absolutely irreducible its image must contain $SL_2(\mathbb F_\ell)$.  Such groups with a complex representation have $\ell \leq 5$.

\begin{proof}  Let $H \subseteq I_\ell$ be the kernel of $\widetilde \rho |_{I_\ell}$ so that the restriction $\rho |_H$ is scalar for some character $\chi: H \to \overline {\mathbb Q}_\ell^\times$ of finite order.  We will choose special extensions of $\chi$ to $G_{\mathbb Q}$.  Say that $\widetilde \rho(G_\ell)$ has order $\ell^n$ i.e. the image of the quotient $G_\ell / H$ has order $\ell^n$, and let $\sigma_0 \in G_\ell$ be a generator. The characteristic polynomial of $\rho(\sigma_0)$ is $( x-\alpha ) \, ( x - \beta )$ where $\alpha/\beta$ is a primitive $\ell^n$th root of unity, and $\rho(\sigma_0)$ is similar to its rational canonical form
\[ \begin{pmatrix} 0 & \alpha \, \beta  \\ -1 & \alpha + \beta \end{pmatrix} = A_\alpha \begin{pmatrix} \beta & 1 \\ & \alpha \end{pmatrix} A_\alpha^{-1} = A_\beta \begin{pmatrix} \alpha & 1 \\ & \beta \end{pmatrix} A_\beta^{-1} \quad \text{where} \quad A_x = \begin{pmatrix} x & -1 \\ 1 & 0 \end{pmatrix}. \]

\noindent (Since $\overline \rho$ is wildly ramified the residual image cannot be diagonal.)  In particular, $\rho(G_\ell)$ may be conjugated to be upper-triangular.  Choose characters $\chi_1, \, \chi_2: G_{\mathbb Q} \to \overline {\mathbb Q}_\ell^\times$ such that 
\begin{enumerate}
\item $\chi_1 |_H = \chi_2 |_H = \chi$;
\item $\chi_1 \left( {\sigma_0}^m \right) = \alpha^m$ and $\chi_2 \left( {\sigma_0}^m \right) = \beta^m$ where $m = \left[ \widetilde \rho(G_\ell): \widetilde \rho(I_\ell) \right]$; and
\item $\chi_1$ and $\chi_2$ have finite order and are unramified outside of $\ell$.
\end{enumerate}

\noindent Now define the characters $\chi_\ell = {\chi_1}^{-1} \chi_2$ and $\chi_0 = \det \rho / (\chi_1 \, \chi_2)$ as well as the representations $\rho_i = {\chi}_i^{-1} \otimes \rho$ for $i = 1, \, 2$.  Then $\chi_\ell$ is wildly ramified at $\ell$, $\chi_0$ is unramified at $\ell$,
\[ \rho_1 |_{I_\ell} \simeq \begin{pmatrix} \chi_\ell & \ast \\ & 1 \end{pmatrix} \quad \text{and} \quad \rho_2 |_{I_\ell} \simeq \begin{pmatrix} {\chi_\ell}^{-1} & \ast \\ & 1 \end{pmatrix}. \]

\noindent Each $\rho_i$ is ordinary and ramified at finitely many primes; while $\overline \rho_i$ is modular, absolutely irreducible, and wildly ramified at $\ell$.  Theorem \ref{wiles} implies that $\rho_i$ and hence $\rho$ are all $\ell$-adically modular.

Let $f(\tau) = \sum_n a_n \, q^n$ denote the $\ell$-adic form associated with $\rho_1$ and $g(\tau) = \sum_n b_n \, q^n$ denote the $\ell$-adic form associated with $\rho_2$.  We have
\begin{enumerate}
\item $f$ and $g$ are ordinary cusp forms of weight $\kappa = 1$;
\item $f$ and $g$ have nebentype $\det \rho_1 = \chi_0 \cdot \chi_\ell$ and $\det \rho_2 = \chi_0 \cdot {\chi_\ell}^{-1}$;
\item $\rho_1 = \chi_\ell \otimes \rho_2$ i.e. $a_p = \chi_\ell \left( \text{Frob}_p \right) \cdot b_p$ for almost all $p \neq \ell$; and
\item $a_\ell = b_\ell = \chi_0 \left( \text{Frob}_\ell \right)$.
\end{enumerate}

\noindent Note that the projective representations $\widetilde \rho \simeq \widetilde \rho_f \simeq \widetilde \rho_g$, so they all have the same prime-to-$\ell$ part of the Artin conductor $N_\emptyset$.  It will follow from \cite[Theorem 11.1]{MR1937198} that such forms are classical modular forms once we verify that $N_\emptyset \geq 5$.  The kernel of $\widetilde \rho: G_{\mathbb Q} \to PGL_2(\mathbb C)$ fixes an extension $L/\mathbb Q$ such that $\text{Gal}(L/\mathbb Q)$ contains $PSL_2(\mathbb F_\ell)$; the only possibilities for the Galois group are $A_4 \simeq PSL_2(\mathbb F_3)$, $S_4 \simeq PGL_2(\mathbb F_3)$, and $A_5 \simeq PSL_2(\mathbb F_5)$. But consulting the tables at \cite{J3} and \cite{J1} we see that there are no such number fields with prime-to-$\ell$ discriminant $N_\emptyset \leq 4$. \end{proof}


\section{Icosahedral Galois Representations}

In this section we specialize to $\ell=5$.  For any continuous complex Galois representation $\rho$, there is a finite extension $L/\mathbb Q$ which makes the following diagram commute:
\[ \begin{CD} 1 @>>> \text{Gal} \left( \overline {\mathbb Q} / L \right) @>>> \text{Gal} \left( \overline {\mathbb Q} / \mathbb Q \right) @>>> \text{Gal} \left( L / \mathbb Q \right) @>>> 1 \\ @. @VVV @VV{\rho}V @VV{\widetilde \rho}V \\ 1 @>>> \mathbb C^{\times} @>>> GL_2 \left( \mathbb C \right) @>>> PGL_2 \left( \mathbb C \right) @>>> 1 \end{CD} \]

\noindent The various possibilities for $\widetilde \rho (G_{\mathbb Q})$ were classified in \cite{MR96g:01046}: it is either cyclic, dihedral, tetrahedral ($A_4 \simeq PSL_2(\mathbb F_3)$), octahedral ($S_4 \simeq PGL_2(\mathbb F_3)$), or icosahedral ($A_5 \simeq PSL_2(\mathbb F_5)$).  We will focus on the latter, so that $L/\mathbb Q$ is an $A_5$-extension.  Our goal in this section is to find a family of Galois representations $\rho$ which satisfy the hypotheses of Theorem \ref{buzzard} by focusing on the extensions $L / \mathbb Q$.  We introduce elliptic curves to give a family residually modular representations.

\subsection{Quintics and Elliptic Curves}

We begin by showing the explicit relationship between a certain class of quintics and elliptic curves.  Most of the exposition that follows in this section is motivated by both \cite{MR96g:01046} and \cite{MR98e:11067}.

The group $A_5$ may be realized as the group of rotations of the icosahedron.  This Platonic Solid has 12 vertices which may be inscribed on the unit sphere, so that after projecting to the extended complex plane we have a natural group action on the complex numbers
\[ \left( \zeta_5 + \zeta_5^4 \right) \, \zeta_5^\nu, \quad \left( \zeta_5^2 + \zeta_5^3 \right) \, \zeta_5^\nu \quad \text{for $\nu \in \mathbb F_5$}; \qquad \infty; \quad \text{and} \quad 0; \]

\noindent generated by the three fractional linear transformations
\[ S \, z = \zeta_5 \, z, \qquad T \, z = \frac {\varepsilon \, z + 1}{z - \varepsilon}, \qquad \text{and} \qquad U \, z = - \frac 1z; \]

\noindent where $\zeta_5$ is a primitive fifth root of unity and $\varepsilon = \zeta_5 + \zeta_5^4$ is a fundamental unit.  Consider the rational functions
\[ \lambda(z) = \frac {\left[ z^2 + 1 \right]^2 \left[ z^2 - 2 \, \varepsilon \, z - 1 \right]^2 \left[ z^2 + 2 \, \varepsilon^{-1} \, z - 1 \right]^2 }{- z \, \left( z^{10} + 11 \, z^5 - 1 \right)}, \quad \mu(z) = \frac {- 125 \, z^5}{z^{10} + 11 \, z^5 - 1}. \]

\noindent The first has nontrivial action by $S$ so we may associate a polynomial of degree 5, while the second has trivial action by $S$ so we may associate a polynomial of degree 6.  The following is the fundamental result in \cite{MR96g:01046}, with the last statement being motivated by the degree 6 resolvent first considered in \cite{MR98e:11067}.

\begin{klein} \label{klein} Fix a quintic $q(x) = x^5 + A \, x^2 + B \, x + C$ over $\mathbb Q$ such that $A^4 - 5 B^3 + 25 \, A \, B \, C$ is nonzero, and denote $L$ as its splitting field.
\begin{enumerate}
\item The rational function
\[ j(z) = (\lambda +3 )^3 \, (\lambda^2 + 11 \, \lambda + 64 ) = \frac {(\mu^2 + 10 \, \mu + 5 )^3}{\mu} \]
\noindent is invariant under $A_5 \simeq \langle S, \, T, \, U \rangle$.  For each $m, \, n \in \overline {\mathbb Q}$, the resolvents
\[ \label{r_substitution} x_{\nu} = \frac m{\lambda(\zeta_5^\nu z) + 3} + \frac n{\left[ \lambda(\zeta_5^\nu z) + 3 \right] \left[ \lambda(\zeta_5^\nu z)^2 + 10 \, \lambda(\zeta_5^\nu z) + 45 \right]} \]
\noindent are roots of the quintic $x^5 + A_{m,n,j} \, x^2 + B_{m,n,j} \, x + C_{m,n,j}$ where
\[ \begin{aligned} A_{m,n,j} & = - \frac {20}j \, \left[ \left( 2 \, m^3 + 3 \, m^2 \, n \right) + 432 \, \frac {6 \, m \, n^2 + n^3}{1728 - j} \right], \\ B_{m,n,j} & = - \frac 5j \, \left[ m^4 -  864 \, \frac {3 \, m^2 \, n^2 + 2 \, m \, n^3}{1728 - j} +  559872 \, \frac {n^4}{(1728 - j)^2} \right], \\ C_{m,n,j} & = - \frac 1j \, \left[ m^5 - 1440 \, \frac {m^3 \, n^2}{1728 - j} + 62208 \, \frac {15 \, m \, n^4 + 4 \, n^5}{(1728 -  j)^2} \right]. \end{aligned} \]
\item There exist $m, \, n, \, j \in \mathbb Q( \sqrt{5 \, \text{Disc}(q)})$ such that $A=A_{m,n,j}$, $B=B_{m,n,j}$ and $C=C_{m,n,j}$.
\item There exists a curve $E$ over $\mathbb Q(\sqrt{5 \, \text{Disc}(q)})$ such that $L(\sqrt{5 \, \text{Disc}(q)})$ is the field generated by sum $x_P+x_{2P}$ of $x$-coordinates of the 5-torsion of $E$.  
\end{enumerate}
\end{klein}

\begin{proof}  The properties stated in the first part are easily verified using a symbolic calculator.  For the second, define $\delta$, $\gamma_4$, $\gamma_6$ and $\text{Disc}(q)$ by
\[ \begin{aligned} 
5^4 \cdot \delta & = A^4 - 5 B^3 + 25 \, A \, B \, C \\ {12}^2 \, 5^5 \cdot \gamma_4 & = 128 \, A^4 \, B^2 - 192 \, A^5 \, C - 600 \, A \, B^3 \, C + 1000 \, A^2 \, B \, C^2 \\ & \quad - 144 \, B^5 + 3125 \, C^4 \\ {12}^3 \, 5^{10} \cdot \gamma_6 & = 1728 \, A^{10} + 10400 \, A^6 \, B^3 + 405000 \, A^2 \, B^6 - 180000 \, A^7 \, B \, C \\ & \quad - 1170000 \, A^3 \, B^4 \, C + 1725000 \, A^4 \, B^2 \, C^2 - 1800000 \, A^5 \, C^3 \\ & \quad + 2812500 \, A \, B^3 \, C^3 - 4687500 \, A^2 \, B \, C^4 \\ & \quad - 2025000 \, B^5 \, C^2 - 9765625 \, C^6 \\ \text{Disc}(q) & = -27 \, A^4 \, B^2 + 108 \, A^5 \, C - 1600 \, A \, B^3 \, C + 2250 \, A^2 \, B \, C^2 \\ & \quad + 256 \, B^5 + 3125  \, C^4 \end{aligned} \]

\noindent By eliminating $m$ and $n$ in the system above, we find that $j$ is a root of the equation
\begin{equation} \label{j-equation} \delta^5 \, j^2 - 1728 \, \left( \gamma_4^3 - \gamma_6^2 + \delta^5 \right) \, j + {1728}^2 \, \gamma_4^3 = 0. \end{equation}

\noindent The statement follows since \eqref{j-equation} has discriminant $5 \cdot \text{Disc}(q)$ and $m$ and $n$ may be expressed in terms of $j$.  For the third statement, fix $j$ as a root of \eqref{j-equation} so that $L(\sqrt{5 \, \text{Disc}(q)})$ is the splitting field of $q'(\mu) = (\mu^2 + 10 \, \mu + 5 )^3 - j \, 
\mu$, and let $E$ be denote the elliptic curve $y^2 = x^3 + 3 \, j / (1728-j) \, x + 2 \, j / (1728-j)$ with invariant $j$.  Given a 5-torsion point $P$ on $E$,
\[ x_P + x_{2P} = -2 \, \frac {\mu^2 + 10 \, \mu + 5}{\mu^2 + 4 \, \mu - 1} \quad \text{for some root of $q'(\mu)$} \]

\noindent so that $L(\sqrt{5 \, \text{Disc}(q)})$ contains the field generated by the sum of the $x$-coordinates of the 5-torsion, and conversely, given a root of $q'(\mu)$,
\[ \mu = \frac {31104 \, x^3}{(x-2)^5 \, j - 1728 \, x^3 \, (x^2 -10 \, x + 34)}  \qquad \text{where $x = x_P + x_{2P}$} \] 

\noindent for some 5-torsion $P$ on $E$.  Hence $L(\sqrt{5 \, \text{Disc}(q)})$ is generated as claimed. \end{proof}

\subsection{$\mathbb Q$-Curves}

We specialize to a family of quintics with particularly nice properties.  The motivation for the following result comes from the classical work of Bring and Jerrard that any quintic can be brought to the form $x^5 + B \, x + C$.  We say an elliptic curve $E$ is associated to $q(x)$ whenever statement (3) of Proposition \ref{klein} holds.  We remind the reader that a $\mathbb Q$-curve is an elliptic curve without complex multiplication which is isogenous to each of its Galois conjugates, as first considered in \cite{MR94g:11042}.  

\begin{q-curve} \label{q-curve} For $t \in \mathbb Q^{\times}$, define the quintic and the curve
\[ q_t(x) = x^5 + 5 \left( \frac {9-5 \, t^2}{5 \, t^2} \right) x + 4 \left( \frac {9-5 \, t^2}{5 \, t^2} \right); \quad E_t: \ y^2 = x^3 + 2 \, x^2 + \frac {3+\sqrt{5} \, t}{2 \sqrt{5} \, t}. \]

\noindent We have the following.
\begin{enumerate}
\item $q_t(x)$ has Galois group contained in $A_5$.  If $t$ is square of a rational number then $q_t(x)$ has Galois group $A_5$.
\item If $t$ is the square of a 5-adic unit, then the decomposition, inertia, and wild inertia groups at 5 are cyclic of order 5.  All higher ramification groups are trivial.
\item$E_t$ is a 2-isogenous $\mathbb Q$-curve that is associated to $q_t(x)$.  Moreover, the isogeny is defined over $\mathbb Q(\sqrt{-2}, \, \sqrt{5})$. 
\item Given a quintic $q(x) = x^5 + B \, x + C$ over $\mathbb Q$ with Galois group $A_5$, the curve $E_t$ is associated to $q(x)$ when $t = 75 \, C^2 / \sqrt{\text{Disc}(q)}$.
\end{enumerate}
\end{q-curve}

Before we show the proof, we show an application of this result. Let $L/\mathbb Q$ be an $A_5$-extension which is unramified outside of $\{2, \, 5, \, \infty \}$.  By \cite{J1}, there are only five such $A_5$-extensions of $\mathbb Q$, so it is a computational exercise to show that the splitting fields actually come from principal quintics: 
\[ \begin{matrix} \text{Original Quintic} & \text{Principal Quintic} & \text{Parameter $t$} \\ \hline x^5 + 20 \, x - 16 & 4 \, x^5 - 25 \, x + 50 & 3/5, \, 15/11 \\ x^5 + 10 \, x^3 - 10 \, x^2 + 35 \, x - 18 & 5 \, x^5 + 20 \, x + 16 & 1 \\ x^5 - 10 \, x^3 + 20 \, x^2 + 110 \, x - 116 & 5 \, x^5 - 20 \, x + 16 & 3 \\ x^5 + 10 \, x^3 - 40 \, x^2 + 60 \, x - 32 & 5 \, x^5 - 5 \, x + 4 & 3/2 \\ x^5 - 10 \, x^3 - 20 \, x^2 + 10 \, x + 216 & 5 \, x^5 + 5 \, x + 8 & 4/3
\end{matrix} \]

\noindent The associated $\mathbb Q$-curves can now be found via Proposition \ref{q-curve}.  In particular, the quintic $x^5 + 10 \, x^3 - 10 \, x^2 + 35 \, x - 18$ first studied in \cite{MR58:22019} has associated $\mathbb Q$-curve
\[ y^2 = x^3 + \left( 5 - \sqrt{5} \right) x^2 + \sqrt{5} \, x. \]

\noindent For more information, see \cite{MR1981898}.  

\begin{proof}  The discriminant $\text{Disc}(q_t) = 2^8 3^2 \, t^{-10} \left( 9- 5 \, t^2 \right)^4$ is a square, so its Galois group is contained in $A_5$.  By \cite{MR83i:12010}, the only solvable trinomials in the form $x^5 + B \, x + C$ with square determinant satisfy
\[ B = 20 \, \frac {(v^2 + v - 1) \, (v^2 - v - 1)}{(v^2 + 1)^2} \, w^4, \quad C = 16 \, \frac {(v^2 + v - 1) \, (2 \, v^2 + 3 \, v - 2)}{(v^2 + 1)^2} \, w^5; \]

\noindent  for some rational $v$ and $w$.  If $q_t(x)$ is solvable by radicals we must have
\[ w = \frac {v^2 - v - 1}{2 \, v^2 + 3 \, v - 2} \quad \text{and} \quad t = \frac {3 \, (v^2 + 1) \, (2 \, v^2 + 3 \, v - 2)^2}{5 \, (2 \, v^3 + 2 \, v^2 - v + 1) \, (v^3 + v^2 + 2 \, v - 2)}. \]

\noindent Say $t = u^2$ for some $u \in \mathbb Q^{\times}$; then the Galois group is properly contained in $A_5$ if and only if we have a rational point on the hyperelliptic curve
\[ y^2 = 15 \, (x^2 + 1) \, (2 \, x^3 + 2 \, x^2 - x + 1) \, (x^3 + x^2 + 2 \, x - 2). \]

\noindent This curve has no rational points as verified using \texttt{MAGMA} \cite{MR1413180} so the quintic has nonsolvable Galois group.

Now we consider the ramification groups at 5.  Fix an embedding $\overline {\mathbb Q} \hookrightarrow \overline {\mathbb Q}_5$, and say $t = u^2$ for some $u \in \mathbb Z_5^{\times}$.  Then
\[ y = \frac {4 \, u}{5 \, \sqrt[4]{ \left( 1 + 5 ( u^4 - 2) \right)}} \in \mathbb Q_5 \implies q_t \left( \frac {x}{5 \, y/4} \right) = \frac {x^5 - x - y}{\left( 5 \, y / 4 \right)^5}; \]

\noindent so that Galois group over $\mathbb Q_5$ is that of the Artin-Schreier quintic $x^5 - x - y$.  As $y$ has 5-adic valuation -1, the results in \cite[page 72, exercise 5]{MR82e:12016} show that the decomposition, inertia, and wild inertia groups at 5 are cyclic of order 5 while the higher ramification groups are trivial.

Set $r = (3+\sqrt{5} \, t) / (2 \sqrt{5} \, t)$ so that $E_t$ is in the form $y^2 = x^3 + 2 \, x^2 + r \, x$.  This curve is singular only if $t = 0, \, \pm 3/\sqrt{5}$ which never happens when $t \in \mathbb Q^{\times}$ so $E_t$ is indeed an elliptic curve.  Denote $\sigma \in \text{Gal} \left( \mathbb Q(\sqrt{-2}, \sqrt{5}) / \mathbb Q(\sqrt{-2}) \right)$ as the nontrivial automorphism, and consider the isogeny $\phi: E_t \to E_t^{\sigma}$ defined by
\[ \phi \left( x, \, y \right) = \left( \frac 1{\sqrt{-2}^2} \, \frac {y^2}{x^2}, \, \frac 1{\sqrt{-2}^3} \, \frac {y \, (r-x^2)}{x^2} \right) \implies \phi \circ \phi^{\sigma} = \left[ -2 \right]. \]

\noindent It follows that $E_t$ is 2-isogenous to its Galois conjugate, so it is a $\mathbb Q$-curve.  The fact that $E_t$ is associated to $q_t(x)$ comes from checking the system of equations in Proposition \ref{klein}.  It is easy to check that $j(E_t)$ is a solution to the quadratic equation in \eqref{j-equation} when $t = 75 \, C^2 / \sqrt{\text{Disc}(q)}$, so $E_t$ is indeed an elliptic curve over $\mathbb Q(\sqrt{5})$ which is associated to $q(x)$.  \end{proof}

\subsection{Geometric Galois Representations}

The extension $L / \mathbb Q$ is associated to the projective representation $\widetilde \rho$ so we construct, under suitable hypotheses, a canonical complex Galois representation as a lift.  We exploit the fact that mod 5 representations of $\mathbb Q$-curves defined over $\mathbb Q(\sqrt{5})$ fit into the composition
\begin{equation} \begin{CD}
    \text{Gal} \left( \overline {\mathbb Q} / \mathbb Q(\sqrt{5}) \right) @>{\overline \rho_{E,5}}>> Z(\mathbb 
    F_5) \cdot SL_2(\mathbb F_5) @>{\pi}>> GL_2(\mathbb C)
\end{CD} \end{equation} 

\noindent where $\pi$ is some finite group representation and $Z(\mathbb F_5)$ is the center of $GL_2(\mathbb F_5)$.  Recall that $\varepsilon = \zeta_5 + \zeta_5^4$ is a fundamental unit in $\mathbb Q(\sqrt{5})$.

\begin{ell_rep} \label{ell_rep}
Let $E = E_t$ be a $\mathbb Q$-curve as in Proposition \ref{q-curve}.  
\begin{enumerate}
\item There exists a faithful representation $\pi: Z(\mathbb F_5) \cdot SL_2(\mathbb F_5) \to GL_2(\mathbb C)$ and a prime ideal $\lambda \subseteq \mathbb Z[\varepsilon, \, \sqrt{-1}]$ above 5 such that $\pi \equiv 1 \mod \lambda$.
\item There exists a 1-dimensional representation $\omega : \text{Gal} \left( \overline {\mathbb Q} / \mathbb Q(\sqrt{5}) \right) \to \mathbb C^{\times}$ such that the twists $\rho_E^{(1)} = \omega \otimes (\pi \circ \overline {\rho}_{E,5})$ and $\rho_E^{(2)} = \omega \otimes \rho_{E,5}$ are continuous Galois representations, ramified at finitely many primes, which are restictions of 2-dimensional representations of $G_{\mathbb Q}$.
\item As representations of $G_{\mathbb Q}$, the determinant $\det \rho_E^{(\kappa)} = \chi \cdot (\varepsilon_5 / \omega_5 )^{\kappa - 1}$ for $\kappa = 1, \, 2$; where $\chi$ is the nontrivial quadratic character which factors through $\text{Gal} \bigl( \mathbb Q(\sqrt{-1}) / \mathbb Q \bigr)$, and $\omega_5$ is a character of order 4 which factors through $\text{Gal} \bigl( \mathbb Q(\zeta_5) / \mathbb Q \bigr)$.  In particular, $\rho_E^{(1)}$ is odd.
\item If $q_t(x)$ has Galois group $A_5$, then the residual representation $\overline \rho_E \equiv \rho_E^{(1)} \equiv \rho_E^{(2)} \mod \lambda$, as a representation $G_{\mathbb Q}$, has image $Z(\mathbb F_5) \cdot SL_2(\mathbb F_5)$.  In particular, the restriction of $\overline \rho_E$ to $\text{Gal} \left( \overline {\mathbb Q} / \mathbb Q(\sqrt{5}) \right)$ is absolutely irreducible.
\end{enumerate}
\end{ell_rep}

A similar result is true for all $\mathbb Q$-curves; see \cite{MR1844206} for proofs using Galois cohomology.  Note that $\rho_E^{(1)}: G_{\mathbb Q} \to GL_2(\mathbb C)$ is a complex representation, whereas $\rho_E^{(2)}: G_{\mathbb Q} \to GL_2(\overline {\mathbb Q}_5)$ is a 5-adic representation.

\begin{proof}  Let $\varpi = (\omega_5(2) \cdot \varepsilon^{-1} - 1) \in \mathbb Q(\sqrt{5}, \, \sqrt{-1})$ in terms of the fundamental unit $\varepsilon$ and the Teichmueller character $\omega_5: \mathbb Z_5^{\times} \to \mathbb C^{\times}$, and $\lambda = \varpi \cdot \mathbb Z[\varepsilon, \, \sqrt{-1}]$ be a prime above 5.  Define the map $\pi: \, Z(\mathbb F_5) \cdot SL_2(\mathbb F_5) \to 
GL_2(\mathbb C)$ by
\[ \pi: \quad \begin{matrix} \begin{pmatrix} 1 & 1 \\ & 1 \end{pmatrix} \\ \\ \begin{pmatrix} & -1 \\ 1 \end{pmatrix} \\ \\ \begin{pmatrix} a & \\ & d \end{pmatrix} \end{matrix} \quad \mapsto \quad \begin{matrix} S & = & \tfrac 12 \begin{pmatrix} \varepsilon & \varpi + 2 \\ \varpi & \varepsilon \end{pmatrix} \\ \\ T & = & \begin{pmatrix} & -1 \\ 1 \end{pmatrix} \\ \\ U(a,d) & = & \begin{pmatrix} \omega_5(a) & \\ & \omega_5(d) \end{pmatrix}. \end{matrix} \]

\noindent We take for granted that $Z(\mathbb F_5) \cdot SL_2(\mathbb F_5)$ is generated by the three matrices on the left, so it suffices to check $\pi$ on these generators.  We follow the exposition in \cite[page 405]{MR97k:11080}.  It is a rather straightforward exercise to show that the orders of each 
generator match correctly i.e. $S^5 = T^4 = U^4 = 1_2$, so it suffices to consider the commutation relations among the image of these generators:
\begin{enumerate}
\item $T \cdot U(a,d) \cdot T^{-1} = U(d,a)$,
\item $U(a,d) \cdot S \cdot U(a,d)^{-1} = S^n$ where $n \equiv a \, d^{-1} \pod{5}$,
\item $T \, S^d \, T^{-1} = U(a,d) \ S^{-d} \, T^{-1} \, S^{-a}$ where $a \, d \equiv 1 \pod{5}$.
\end{enumerate}

\noindent  (Condition 2 is {\em not} satisfied if $n \equiv \pm 2 \pod{5}$; this explains why we must restrict to $Z(\mathbb F_5) \cdot SL_2(\mathbb F_5)$ where $n \equiv \pm 1$.)  It is also straightforward to verify these relations.  The the congruence $\pi \equiv 1 \pod \lambda$ follows because $2 - \varepsilon = \varepsilon^2 \, \varpi \, \varpi^c$, $2 - \omega_5(2) = \varepsilon  \, \varpi \, (\varepsilon \, \varpi^c - 1)$, and $\sqrt{5} = \varepsilon \, \varpi \, \varpi^c$ are elements of $\lambda$.  (Recall $c$ is complex conjugation.)

We now prove statement (2) of the proposition.  By Proposition \ref{q-curve}, the isogeny for $E$ is defined over $\mathbb Q(\sqrt{-2}, \sqrt{5})$ so $a_{\sigma(\mathfrak p)} = \left( -2/\mathbb N \, \mathfrak p \right) \cdot a_{\mathfrak p}$ in terms of the trace of Frobenius of $\rho_{E,5}$, the Legendre symbol $(-2/\ast)$, the norm map $\mathbb N$, and the nontrivial automorphism $\sigma \in \text{Gal} \left( \mathbb Q(\sqrt{5}) / \mathbb Q \right)$.  Hence the twisted representations $\rho_E^{(\kappa)}$ will 
be Galois invariant and defined over $\mathbb Q$ if we exhibit a character $\omega: \mathbb Z[\varepsilon] \to \mathbb C$ such that $\omega^{\sigma - 1} = (-2 / \ast) \circ \mathbb N$.  The field  $\mathbb Q(\sqrt{5})$ has narrow class number 1 so we extend the Dirichlet characters of conductors $4 \cdot \mathbb Z[\varepsilon]$, $8 \cdot \mathbb Z[\varepsilon]$, and $\sqrt{5} \cdot \mathbb Z[\varepsilon]$ as follows:
\[ \omega_4 : \left \{ \begin{matrix} -1 \\ \varepsilon \end{matrix}  \ \mapsto \ \begin{matrix} -1 \\ \zeta_6 \end{matrix} \right \}, \qquad \omega_8 : \left \{ \begin{matrix} -1 \\ 1+4 \, \varepsilon \\ \varepsilon \end{matrix} \ \mapsto \ \begin{matrix} -1 \\ -1 \\ \zeta_{12} \end{matrix} \right \}, \qquad \omega_5 : \varepsilon \mapsto \zeta_4. \]

\noindent It is straightforward to verify on the generators that
\[ \omega = {\omega_4}^3 \, {\omega_8}^3 \, \omega_5 \implies \omega^{\sigma - 1} = \left( \frac {-2 \ }{\ast} \right) \circ \mathbb N \quad \text{and} \quad \omega^2 = \left[ \left( \frac {-1 \ }{\ast} \right) \cdot {\omega_5}^{-1} \right] \circ \mathbb N. \]

\noindent Since $\omega$ is trivial on the totally positive units $\varepsilon^{2 \, n}$, we identify $\omega$ with a Galois representation through the composition
\[ \begin{CD} \text{Gal} \left( \overline {\mathbb Q} / \mathbb Q(\sqrt{5}) \right) @>>> \left( \mathbb Z[\varepsilon]  \left / 8 \sqrt{5} \cdot \mathbb Z[\varepsilon] \right. \right)^{\times} @>{\omega}>> \mathbb C^{\times} \end{CD} \]

\noindent via Class Field Theory.  It is clear that $\rho_E^{(\kappa)}$ is continuous and ramified at finitely many primes because the same is true for both $\rho_{E,5}$ and $\omega$.

We have $\omega^2 = \chi \cdot \omega_5^{-1}$ with $\chi$ the nontrivial quadratic character ramified at 2.  It is well-known that $\det \rho_{E,5} = \varepsilon_5$ so that
\[ \det \rho_E^{(1)} = \omega^2 \cdot  \det \pi \circ \overline \rho_{E,5} = \chi \cdot ( \omega_5 / \omega_5), \quad \det \rho_E^{(2)} = \omega^2 \cdot \det \rho_{E,5} = \chi \cdot ( \varepsilon_5 / \omega_5); \]

\noindent hence statement (3) follows.

To show statement (4) note that $Z(\mathbb F_5) \cdot SL_2(\mathbb F_5)$ is the subgroup of $GL_2(\mathbb F_5)$ of index 2 consisting of those matrices with square determinants.  Since $\det \overline \rho_E \equiv \chi \pod 5$ is a quadratic character, the image of $\overline \rho_E$ is contained in $Z(\mathbb F_5) \cdot SL_2(\mathbb F_5)$.  If $q_t(x)$ has Galois group $A_5$, the restriction $\overline \rho_E  |_{\mathbb Q(\sqrt{5})}$ is the twist of an absolutely irreducible mod 5 representation, which has image $Z(\mathbb F_5) \cdot SL_2(\mathbb F_5)$.  This gives equality. \end{proof}

\subsection{Modularity Results}

We now apply the results from the previous section.  Note that the following result imposes no local conditions at 3.

\begin{modularity} \label{modularity} Let $E = E_t$ be a $\mathbb Q$-curve as in Proposition \ref{q-curve}.  Then $E$ is modular.  In particular, both $\rho_E^{(2)}$ and $\overline \rho_E$ are modular.
\end{modularity}

\begin{proof}  We use the arguments in \cite{MR1639612} and \cite{MR1844206} on the modularity of representations which are potentially Barsotti-Tate, so we consider the local properties of $E$ at $\ell = 3$.  

For a prime $\ell$, denote $\rho_\ell: G_{\mathbb Q} \to GL_2(\overline {\mathbb Q}_\ell)$ as the $\ell$-adic representation associated to the Galois invariant twist $\omega \otimes \rho_{E,\ell}$ as in Proposition \ref{ell_rep}.  When $\ell = 5$ we identify $\rho_\ell = \rho_E^{(2)}$, but it suffices to show $\rho_\ell$ is modular for $\ell = 3$.  Recall by \cite[Proposition 2.10]{MR1844206} that there exists an abelian variety $A$ defined over $\mathbb Q$ such that there exists an embedding $\mathbb Z[\sqrt{-2}] \hookrightarrow \text{End}(A) \otimes_{\mathbb Z} \mathbb Q$ and $\rho_3 \simeq \rho_{A,\lambda}$ can be realized as the $\lambda$-adic representation of the Tate module
\[  \projlim_{n \to \infty} A \left[ \lambda^n \right] \otimes_{\mathbb Z_\lambda} \mathbb Q_\lambda \qquad \text{where} \qquad \lambda = \left( 1 + \sqrt{-2} \right) \cdot \mathbb Z[ \sqrt{-2} ]. \]

\noindent The elliptic curve $E$ has a Weierstrass equation of the form $y^2 = x^3 + 2 \, x^2 + r \, x$ where $r = (3+\sqrt{5} \, t)/(2 \sqrt{5} \, t)$ for some $t \in \mathbb Q^{\times}$.  Fix an embedding $\overline {\mathbb Q} \hookrightarrow \overline {\mathbb Q}_3$ as well as a normalized valuation $\nu: \mathbb Q_3 \to \mathbb Z$.  If $\nu(t) \leq 1$ then both $r, 1-r \in \mathbb Z_3[\sqrt{5}]^\times$ so that $E$ has good reduction at the primes in $\mathbb Q \bigl( \sqrt{-2}, \sqrt{5} \bigr)$ above 3; \cite[Theorem 5.2]{MR1844206} shows $E$ is modular.  (Modularity essentially follows from \cite[Theorem 5.3]{MR97d:11172} because $\rho_\ell$ is crystalline.)  Assume then that $\nu(t) \geq 2$, and write $t = 3^{n+2} \, u$ for some $u \in \mathbb Z_3^{\times}$ so that the twist of $E$ by $\sqrt{-3}^{n+1}$ has good reduction at the primes above 3.  If $n \equiv 3 \pod{4}$ then $E$ itself has good reduction at the primes above 3, so $E$ is again modular.  We assume that $n \not \equiv 3 \pod{4}$, and follow the arguments in \cite[\S 3]{MR1844206} to show $E$ is modular.  Note that reduction of the twisted elliptic curve modulo the primes above 3 is $y^2 = x^3 + D \, x$ where $D = (-1)^n/(\sqrt{5} \, u) \in \mathbb F_9^{\times}$, so the variety $A$ has potential supersingular reduction.

Denote $F = \mathbb Q_3(\sqrt{-3})$ as the ramified quadratic extension of $\mathbb Q_3$, and $\psi: G_F = \text{Gal} \left( \overline {\mathbb Q}_3 / \mathbb Q_3(\sqrt{-3}) \right) \to \{ \pm 1 \}$ as the ramified quadratic character; the twisted curve $A^\psi$ defined over $F$ has good supersingular reduction and so the twisted representation $\rho_3 \vert_{G_F} \otimes \psi$ is unramified.  (Compare with \cite[Lemma 3.4]{MR1844206}.)  Clearly the restriction of $\overline \rho_3$ to $\text{Gal} \left( \overline {\mathbb Q} / \mathbb Q(\sqrt{-3}) \right)$ is absolutely irreducible, and so by \cite[Lemma 3.2]{MR1844206} the residual representation $\overline \rho_3$ is modular.  By \cite[Lemma 3.5]{MR1844206} the centralizer of the image $\overline \rho_3 \bigl( G_3 \bigr)$ consists of scalars.  Following the notation in \cite[Appendix B]{MR1639612}, it suffices to show that the representation
\[ \begin{CD} WD(\rho_3): \quad W_{\mathbb Q} \to GL_2(\overline {\mathbb Q}_3)\end{CD} \]

\noindent of the Weil group of $\mathbb Q$ is ``strongly acceptable'' for $\overline {\rho}_3$. (Compare with \cite[Lemma 3.6]{MR1844206}.)  Denote the restriction $\tau = WD(\rho_3) \vert_{I_3}$.  Because twisted representation $\rho_3 \vert_{G_F} \otimes \psi$ is unramified and $\psi$ is a quadratic character we have the identity
\[ WD \left( \rho_3 \vert_{G_F} \otimes \psi \right) \simeq WD \left( \rho_3 \vert_{G_F} \right) \otimes_{\overline {\mathbb Q}_3} WD(\psi) \implies \tau \vert_{I_F} \simeq WD(\psi) \vert_{I_F}. \]

\noindent That is, $\tau \vert_{I_F}$ is a nontrivial quadratic character for $I_F$.  As in Proposition \ref{ell_rep}, the determinant $\det \rho_\ell = \chi \cdot \varepsilon_\ell / \omega_5$ so by \cite[\S B.2]{MR1639612} we have
\[ WD \left( \det \rho_3 \right) \simeq \left( \left( \chi / \omega_5 \right) \vert_{W_F} \otimes_{\mathbb Q_3(\sqrt{-1})} \overline {\mathbb Q}_3 \right) \otimes_{\overline {\mathbb Q}_3} WD(\varepsilon_3). \]

\noindent This character is unramified and $\psi$ is a quadratic character so $\det \tau = 1$.  We conclude that $\tau = \tilde \phi^2 \oplus \tilde \phi^6$ in terms of the Teichmuller lift of the fundamental character $\phi: G_{\mathbb Q} \to \mathbb F_9^\times$ of level 2.  If $F'$ is an extension of of $\mathbb Q_3$ such that $\tau \vert_{I_{F'}}$ is trivial, then $F'$ contains the splitting field of $x^4 - 3$ and so $F \subseteq F'$, so that $\rho_3 \vert_{I_{F'}} \simeq \rho_{A^\psi, \lambda} \vert_{I_{F'}}$ is Barsotti-Tate.   It follows from the criteria of \cite[\S 1.2]{MR1639612} that $\rho_3$ is a deformation of $\overline {\rho}_3$ of type $\tau$, and so it follows from \cite[Corollary 2.3.2]{MR1639612} that $\tau$ is ``acceptable'' for $\overline {\rho}_3$.  Since $A$ is potentially supersingular at 3, we have $\overline {\rho}_3 \vert_{I_3} \otimes \overline {\mathbb F}_3 \simeq \phi \oplus \phi^3$ and so using the criteria in \cite[\S 1.2]{MR1639612} we conclude that $\tau$ is indeed ``strongly acceptable'' for $\overline \rho_3$.  \end{proof}

\begin{goins} 
Let $\widetilde \rho: G_{\mathbb Q} \to PGL_2(\mathbb C)$ be a projective icosahedral Galois representation, and assume that the field $L$ fixed by its kernel is the splitting field of a quintic $x^5 + B \, x + C$ with $75 \, C^2 / \sqrt{256 \, B^5 + 3125 \, C^4}$ the square of a 5-adic unit.  Then for any continuous lift $\rho : G_{\mathbb Q} \to GL_2(\mathbb C)$ of $\widetilde \rho$,

\begin{enumerate}
\item $\rho$ is ramified at finitely many primes;
\item $\overline \rho$ is absolutely irreducible when restricted to $\text{Gal} \left( \overline {\mathbb Q} / \mathbb Q(\sqrt{5}) \right)$, modular, and wildly ramified at 5;
\item $\rho(G_5)$ is finite and $\widetilde \rho(G_5)$ is cyclic of order 5.
\end{enumerate}

\noindent In particular, $\rho$ is (classically) modular.
\end{goins}

By Propositions \ref{q-curve} and \ref{modularity} there are infinitely many projective representations that satisfy the hypotheses of the Theorem; for example, the family of quintics
\[ x^5 + 5 \left( \frac {9-5 \, u^4}{5 \, u^4} \right) x + 4 \left( \frac {9-5 \, u^4}{5 \, u^4} \right), \qquad u \in \mathbb Q \cap \mathbb Z_5^{\times}; \]

\noindent would suffice to yield such representations.

\begin{proof} The modularity of $\rho$ would follow from Theorem \ref{buzzard} if we show the three statements hold true.  Choose a 2-isogenous $\mathbb Q$-curve $E$ associated to $x^5 + B \, x + C$ via Proposition \ref{q-curve}, and consider the representation $\rho_E^{(1)}: G_{\mathbb Q} \to GL_2(\mathbb C)$ as in Proposition \ref{ell_rep}.  We will show that $\rho$ is a twist of $\rho_E^{(1)}$. 

Denote $L^{(1)}$ as the field fixed by the kernel of $\widetilde \rho_E^{(1)}$.  First we show $L^{(1)} = L$.  It is well-known that the field fixed by the kernel of $\overline \rho_{E,5}$ is contained in the field generated by the coordinates of the 5-torsion of $E$, so the field fixed by the kernel of the composition
\[ \begin{CD} \text{Gal} \left( \overline {\mathbb Q} / \mathbb Q(\sqrt{5}) \right) @>{\overline \rho_{E,5}}>> Z(\mathbb F_5) \cdot SL_2( \mathbb F_5) @>>> PSL_2(\mathbb F_5) \end{CD} \]

\noindent must be generated by those coordinates of the 5-torsion which are fixed by action of the scalars $Z(\mathbb F_5)$ i.e. the field generated by sum $x_P+x_{2P}$ of $x$-coordinates of the 5-torsion of $E$.  It follows from Proposition \ref{klein} that this field is $L(\sqrt{5})$.  On the other hand, it follows from Proposition \ref{ell_rep} that $\overline \rho_{E,5}$ is surjective while $\pi$ is injective so $L(\sqrt{5})$ is also the field fixed by the kernel of the composition
\[ \begin{CD} \text{Gal} \left( \overline {\mathbb Q} / \mathbb Q(\sqrt{5}) \right) @>{\rho_E^{(1)}}>> GL_2( \mathbb C) @>>> PGL_2(\mathbb C); \end{CD} \]

\noindent hence $L^{(1)}(\sqrt{5}) = L(\sqrt{5})$.  Note $\sqrt{5} \notin L$ because $L / \mathbb Q$ is a nonsolvable extension and has no subfields of degree 2, so $L^{(1)} \supseteq L$.  Similarly, $L^{(1)} \subseteq L$. Now consider the following diagram:
\[ \begin{CD} 1 @>>> \text{Gal} \left( \overline {\mathbb Q} / L \right) @>>> \text{Gal} \left( \overline {\mathbb Q} / \mathbb Q \right) @>>> \text{Gal} \left( L / \mathbb Q \right) @>>> 1 \\  @. @VV{\chi_E}V @VV{\rho, \ \rho_E^{(1)}}V @VV{\widetilde \rho, \ \widetilde \rho_E^{(1)} }V @. \\  1 @>>> \mathbb C^{\times} @>>> GL_2(\mathbb C) @>>> PGL_2(\mathbb C) @>>> 1 \\
\end{CD} \]

\noindent There are only two such projective representations so upon choosing a Galois conjugate of $\pi$ as in Proposition \ref{ell_rep} we have $\widetilde \rho \simeq \widetilde \rho_E^{(1)}$; without loss of generality we assume equality.  To show $\rho$ is a twist of $\rho_E^{(1)}$ we consider the character $\chi_E$ in the left most column.  For automorphisms $\sigma$ and $\tau$, define the symbol $\chi_E(\sigma)$ by the relation $\chi_E(\sigma) \cdot 1_2 = \rho(\sigma) \cdot \rho_E^{(1)}(\sigma)^{-1}$.  Then $\chi_E$ is actually a multiplicative character since
\[ \begin{aligned} \chi_E(\sigma \, \tau) \cdot 1_2 & = \rho(\sigma \, \tau) \, \rho_E^{(1)}(\sigma \, \tau)^{-1} = \rho(\sigma) \left[ \rho(\tau) \, \rho_E^{(1)}(\tau)^{-1} \right] \rho_E^{(1)}(\sigma)^{-1} \\ & = \chi_E(\sigma) \, \chi_E(\tau) \cdot 1_2. \end{aligned} \]

\noindent  Hence, $\rho \simeq \chi_E \otimes \rho_E^{(1)}$.  This character is of finite order because $\rho$ is continuous i.e. has finite image.  The statements that $\overline \rho$ is wildly ramified at 5 and $\widetilde \rho(G_5)$ is cyclic of order 5 now follow from Proposition \ref{q-curve}, while the statement that $\overline \rho$ is modular follows from Proposition \ref{modularity}.    \end{proof}


\bibliographystyle{plain}

\begin{thebibliography}{10}

\bibitem{MR42:7460}
{\em S\'eminaire {B}ourbaki. {V}ol. 1968/69: {E}xpos\'es 347--363}.
\newblock Springer-Verlag, Berlin, 1971.

\bibitem{MR58:22019}
Joe~P. Buhler.
\newblock {\em Icosahedral {G}alois representations}.
\newblock Springer-Verlag, Berlin-New York, 1978.
\newblock Lecture Notes in Mathematics, Vol. 654.

\bibitem{MR97k:11080}
Daniel Bump.
\newblock {\em Automorphic forms and representations}.
\newblock Cambridge University Press, Cambridge, 1997.

\bibitem{MR1937198}
Kevin Buzzard.
\newblock Analytic continuation of overconvergent eigenforms.
\newblock {\em J. Amer. Math. Soc.}, 16(1):29--55 (electronic), 2003.

\bibitem{MR1845181}
Kevin Buzzard, Mark Dickinson, Nick Shepherd-Barron, and Richard Taylor.
\newblock On icosahedral {A}rtin representations.
\newblock {\em Duke Math. J.}, 109(2):283--318, 2001.

\bibitem{MR1897901}
Kevin Buzzard and William~A. Stein.
\newblock A mod five approach to modularity of icosahedral {G}alois
  representations.
\newblock {\em Pacific J. Math.}, 203(2):265--282, 2002.

\bibitem{MR2000j:11062}
Kevin Buzzard and Richard Taylor.
\newblock Companion forms and weight one forms.
\newblock {\em Ann. of Math. (2)}, 149(3):905--919, 1999.

\bibitem{MR1413180}
John Cannon and Catherine Playoust.
\newblock M{AGMA}: a new computer algebra system.
\newblock {\em Euromath Bull.}, 2(1):113--144, 1996.

\bibitem{MR93i:11063}
Robert~F. Coleman and Jos{\'e}~Felipe Voloch.
\newblock Companion forms and {K}odaira-{S}pencer theory.
\newblock {\em Invent. Math.}, 110(2):263--281, 1992.

\bibitem{MR1639612}
Brian Conrad, Fred Diamond, and Richard Taylor.
\newblock Modularity of certain potentially {B}arsotti-{T}ate {G}alois
  representations.
\newblock {\em J. Amer. Math. Soc.}, 1999.

\bibitem{MR99d:11067b}
Henri Darmon, Fred Diamond, and Richard Taylor.
\newblock Fermat's last theorem.
\newblock In {\em Elliptic curves, modular forms \& Fermat's last theorem (Hong
  Kong, 1993)}, pages 2--140. Internat. Press, Cambridge, MA, 1997.

\bibitem{MR52:284}
Pierre Deligne and Jean-Pierre Serre.
\newblock Formes modulaires de poids $1$.
\newblock {\em Ann. Sci. \'Ecole Norm. Sup. (4)}, 7:507--530 (1975), 1974.

\bibitem{MR97d:11172}
Fred Diamond.
\newblock On deformation rings and {H}ecke rings.
\newblock {\em Ann. of Math. (2)}, 144(1):137--166, 1996.

\bibitem{MR1638490}
Fred Diamond.
\newblock An extension of {W}iles' results.
\newblock In {\em Modular forms and Fermat's last theorem (Boston, MA, 1995)},
  pages 475--489. Springer, New York, 1997.

\bibitem{MR98c:11047}
Fred Diamond.
\newblock The {T}aylor-{W}iles construction and multiplicity one.
\newblock {\em Invent. Math.}, 128(2):379--391, 1997.

\bibitem{MR2002m:11042}
Fred Diamond, Matthias Flach, and Li~Guo.
\newblock The {B}loch-{K}ato conjecture for adjoint motives of modular forms.
\newblock {\em Math. Res. Lett.}, 8(4):437--442, 2001.

\bibitem{MR93h:11124}
Bas Edixhoven.
\newblock The weight in {S}erre's conjectures on modular forms.
\newblock {\em Invent. Math.}, 109(3):563--594, 1992.

\bibitem{MR1844206}
Jordan~S. Ellenberg and Chris Skinner.
\newblock On the modularity of $\mathbb {Q}$-curves.
\newblock {\em Duke Math. J.}, 109(1):97--122, 2001.

\bibitem{MR1981898}
Edray~Herber Goins.
\newblock Icosahedral {$\mathbb Q$}-curve extensions.
\newblock {\em Math. Res. Lett.}, 10(2-3):205--217, 2003.

\bibitem{MR91i:11060}
Benedict~H. Gross.
\newblock A tameness criterion for {G}alois representations associated to
  modular forms (mod $p$).
\newblock {\em Duke Math. J.}, 61(2):445--517, 1990.

\bibitem{MR87k:11049}
Haruzo Hida.
\newblock Galois representations into ${\rm {G}{L}}\sb 2({\bf {Z}}\sb p[[{X}]])$ attached to ordinary cusp forms.
\newblock {\em Invent. Math.}, 85(3):545--613, 1986.

\bibitem{MR88i:11023}
Haruzo Hida.
\newblock Iwasawa modules attached to congruences of cusp forms.
\newblock {\em Ann. Sci. \'Ecole Norm. Sup. (4)}, 19(2):231--273, 1986.

\bibitem{MR94j:11044}
Haruzo Hida.
\newblock {\em Elementary theory of ${L}$-functions and {E}isenstein series}.
\newblock Cambridge University Press, Cambridge, 1993.

\bibitem{MR2002c:11057}
Arnaud Jehanne and Michael M{\"u}ller.
\newblock Modularity of an odd icosahedral representation.
\newblock {\em J. Th\'eor. Nombres Bordeaux}, 12(2):475--482, 2000.
\newblock Colloque International de Th\'eorie des Nombres (Talence, 1999).

\bibitem{J1}
John Jones.
\newblock Cumulative quintic tables.
\newblock {\em \newline{\tt
  http://hobbes.la.asu.edu/Number\_Fields/cum-quintics.html}}, 1999.

\bibitem{J3}
John Jones.
\newblock Cumulative tables of quartic fields.
\newblock {\em \newline{\tt
  http://hobbes.la.asu.edu/Number\_Fields/cum-quartics.html}}, 1999.

\bibitem{MR96a:11128}
Ian Kiming and Xiang~Dong Wang.
\newblock Examples of $2$-dimensional, odd {G}alois representations of ${A}\sb
  5$-type over ${\bf {q}}$ satisfying the {A}rtin conjecture.
\newblock In {\em On Artin's conjecture for odd $2$-dimensional
  representations}, pages 109--121. Springer, Berlin, 1994.

\bibitem{MR18:329c}
Felix Klein.
\newblock {\em Lectures on the icosahedron and the solution of equations of the
  fifth degree}.
\newblock Dover Publications Inc., New York, N.Y., revised edition, 1956.
\newblock Translated into English by George Gavin Morrice.

\bibitem{MR96g:01046}
Felix Klein.
\newblock {\em Vorlesungen \"uber das {I}kosaeder und die {A}ufl\"osung der
  {G}leichungen vom f\"unften {G}rade}.
\newblock Birkh\"auser Verlag, Basel, 1993.
\newblock Reprint of the 1884 original, Edited, with an introduction and
  commentary by Peter Slodowy.

\bibitem{MR98e:11067}
Annette Klute.
\newblock Icosahedral {G}alois extensions and elliptic curves.
\newblock {\em Manuscripta Math.}, 93(3):301--324, 1997.

\bibitem{MR90k:11057}
B.~Mazur.
\newblock Deforming {G}alois representations.
\newblock In {\em Galois groups over ${\bf Q}$ (Berkeley, CA, 1987)}, pages
  385--437. Springer, New York, 1989.

\bibitem{MR1638481}
Barry Mazur.
\newblock An introduction to the deformation theory of {G}alois
  representations.
\newblock In {\em Modular forms and Fermat's last theorem (Boston, MA, 1995)},
  pages 243--311. Springer, New York, 1997.

\bibitem{MR88e:14028}
J.~S. Milne.
\newblock {\em Arithmetic duality theorems}.
\newblock Academic Press Inc., Boston, MA, 1986.

\bibitem{MR91g:11066}
K.~A. Ribet.
\newblock On modular representations of ${\rm {G}al}(\overline{\bf {Q}}/{\bf {Q}})$ arising from modular forms.
\newblock {\em Invent. Math.}, 100(2):431--476, 1990.

\bibitem{MR94g:11042}
Kenneth~A. Ribet.
\newblock Abelian varieties over ${\bf {Q}}$ and modular forms.
\newblock In {\em Algebra and topology 1992 (Taej\u on)}, pages 53--79. Korea
  Adv. Inst. Sci. Tech., Taej\u on, 1992.

\bibitem{MR95d:11056}
Kenneth~A. Ribet.
\newblock Report on mod $\ell$ representations of ${\rm {G}al}(\overline{\bf {Q}}/{\bf {Q}})$.
\newblock In {\em Motives (Seattle, WA, 1991)}, pages 639--676. Amer. Math.
  Soc., Providence, RI, 1994.

\bibitem{MR83i:12010}
G.~Roland, N.~Yui, and D.~Zagier.
\newblock A parametric family of quintic polynomials with {G}alois group
  ${D}\sb{5}$.
\newblock {\em J. Number Theory}, 15(1):137--142, 1982.

\bibitem{MR57:16310}
Jean-Pierre Serre.
\newblock Repr\'esentations $\ell$-adiques.
\newblock In {\em Algebraic number theory (Kyoto Internat. Sympos., Res. Inst.
  Math. Sci., Univ. Kyoto, Kyoto, 1976)}, pages 177--193. Japan Soc. Promotion
  Sci., Tokyo, 1977.

\bibitem{MR82e:12016}
Jean-Pierre Serre.
\newblock {\em Local fields}.
\newblock Springer-Verlag, New York-Berlin, 1979.
\newblock Translated from the French by Marvin Jay Greenberg.

\bibitem{MR47:3318}
Goro Shimura.
\newblock {\em Introduction to the arithmetic theory of automorphic functions}.
\newblock Publications of the Mathematical Society of Japan, No. 11. Iwanami
  Shoten, Publishers, Tokyo, 1971.
\newblock Kano Memorial Lectures, No. 1.

\bibitem{MR1638477}
Lawrence~C. Washington.
\newblock Galois cohomology.
\newblock In {\em Modular forms and Fermat's last theorem (Boston, MA, 1995)},
  pages 101--120. Springer, New York, 1997.

\bibitem{MR96d:11071}
Andrew Wiles.
\newblock Modular elliptic curves and {F}ermat's last theorem.
\newblock {\em Ann. of Math. (2)}, 141(3):443--551, 1995.

\end{thebibliography}

\end{document}